\documentclass[12pt,a4paper,twoside,final,notitlepage, leqno]{article}
\usepackage[english]{babel}
\usepackage[T1]{fontenc}  
\usepackage{epsfig, graphicx, amssymb}
\usepackage{amsmath,amsthm,epsfig,amsfonts}  
\usepackage{float}
\usepackage{color}
\newcommand{\monitem}{ \smallskip \noindent $\bullet$ \quad  } 
\newcommand{\moneq}{\vspace*{-7pt} \begin{equation} \displaystyle } 
\newcommand{\moneqstar}{\vspace*{-6pt} \begin{equation*} \displaystyle } 
\newcommand{\monendstar}{\vspace*{-6pt} \end{equation*}   }
\newcommand{\monend}{\vspace*{-7pt} \end{equation}   }
\newcommand{\moneqarraystar}{ \begin{eqnarray*} \displaystyle } 
\newcommand{\monendarraystar}{ \end{eqnarray*}   }

\newcommand{\dd}{{\rm d}}

\newcommand{\RR}[0]{\mathbb{R}}

\newcommand{\NN}[0]{\mathbb{N}}

\definecolor{vertfonce}{rgb}{0.0, 0.5, 0.0}

\def\section*#1{}

\usepackage{fancyhdr}
\renewcommand{\headrulewidth}{0pt}

\begin{document} 

\fancypagestyle{plain}{ \fancyfoot{} \renewcommand{\footrulewidth}{0pt}}
\fancypagestyle{plain}{ \fancyhead{} \renewcommand{\headrulewidth}{0pt}}

~

  \vskip 2.1 cm

\centerline {\bf \LARGE  Nonlinear fourth order Taylor expansion }

\bigskip
 
\centerline {\bf \LARGE of lattice Boltzmann schemes  }

 \bigskip  \bigskip \bigskip

\centerline { \large    Fran\c{c}ois Dubois$^{ab}$}

\smallskip  \bigskip 

\centerline { \it  \small   
  $^a$   Laboratoire de Math\'ematiques d'Orsay, Facult\'e des Sciences d'Orsay,}

\centerline { \it  \small   Universit\'e Paris-Saclay, France.} 

\centerline { \it  \small   
$^b$    Conservatoire National des Arts et M\'etiers, LMSSC laboratory,  Paris, France.}


\bigskip  \bigskip  

\centerline { 19 June 2020  
  {\footnote {\rm  \small $\,${\it Asymptotic Analysis}, volume 127, pages 297-337, 2022. \\
This contribution has been first  presented
at the Institut Henri Poincar\'e (Paris, France) the 10~January 2018. Edition 22 August 2024.  }}}

 \bigskip \bigskip  
 {\bf Keywords}: partial differential equations, asymptotic analysis  

 {\bf AMS classification}:
 65Q05,   
 76N15,  
 82C20.   
 
 {\bf PACS numbers}:  
02.70.Ns, 
05.20.Dd, 
47.10.+g 

\bigskip  \bigskip   
\noindent {\bf \large Abstract } 

\noindent
We propose a formal expansion of multiple relaxation times lattice Boltzmann schemes in
terms of a single infinitesimal numerical variable.
The result is a system of partial differential equations for the conserved moments
of the lattice Boltzmann scheme.
The expansion is presented  in the nonlinear case up to fourth order accuracy.
The asymptotic corrections of  the nonconserved moments  are developed in terms of  
equilibrium values and  partial differentials of the conserved moments.
Both expansions are coupled and conduct to explicit compact formulas.
The new algebraic expressions are validated with previous results obtained with this framework. 
The example of isothermal D2Q9 lattice Boltzmann scheme illustrates the theoretical framework.

\noindent

\newpage 
%

\bigskip \bigskip    \noindent {\bf \large    1) \quad  Introduction  }    

\fancyhead[EC]{\sc{Fran\c{c}ois Dubois}} 
\fancyhead[OC]{\sc{Nonlinear fourth order Taylor expansion of Lattice Boltzmann schemes}}
\fancyfoot[C]{\oldstylenums{\thepage}}

\noindent 
The lattice Boltzmann schemes in their modern  form have been developed 
with the contributions of d'Humi\`eres, Lallemand, Succi    \cite{HSB91, QHL92, DDH92, LL00}
and many others. 
An underlying lattice Boltzmann equation is discretised on a cartesian grid 
with a finite set of velocities chosen in such a way that during one time step, an exact transport is done 
between two  vertices of the mesh.  
A lattice Boltzmann scheme is composed with two steps:
a nonlinear  local relaxation step,   followed by a linear advection scheme coupling
a given vertex with a given family of neighbours.
The relaxation step follows in general the ``BGK'' approximation 
introduced by Bhatnagar, Gross and Krook \cite{BGK54}.
This numerical method allows the simulation of an important number of physical phenomena as
isothermal flows, compressible flows with heat transfer,  non-ideal fluids, multiphase and multi-component flows, 
microscale gas flows, soft-matter flows,...  up to quantum mechanics. 
For the actual status of the method and the various applications, we refer {\it e.g.} to the books of Guo and Shu  \cite{GS13}, 
Kr\"uger {\it et al.} \cite{K3S2V17} and to the prospective article  of Succi \cite{Su15}. 
The lattice Boltzmann method is founded on a mesoscopic Boltzmann model 
but is not able  in general to solve the Boltzmann equation or associated kinetic models.
It is admitted that essentially macroscopic models are approximated with the lattice Boltzmann
schemes.

\smallskip \noindent
The link between   mesoscopic and  macroscopic models is never straightforward and an analysis is necessary.
The  usual  method of formal analysis is founded on the Chapman-Enskog \cite{CC39} 
expansion  in the spirit of the book of Huang \cite{Hu63}.
The classical approach  for the Boltzmann equation  and other mesoscopic models posed in a continuum space-time.
The  Chapman-Enskog method has been adapted to  lattice Boltzmann schemes with discrete space and  time equations by 
Chen and  Doolen \cite{CD98} and by  Qian and Zhou \cite{QZ2k}.   
It allows the derivation  the equivalent partial differential equations founded on a multiscale time analysis. 
This  classical Chapman Enskog expansion is popular in the
lattive Boltzmann community, as reported in the contributions of  Lallemand and Luo \cite{LL00},
Philippi and Hegele \cite{PH06},
Shan {\it et al.} \cite{SHC06} and Otomo  {\it et al.} \cite{OBD17}
among others.     
A main drawback of this approach is the fact that multiscale expansions
are used without a clear mathematical signification  of the various variables and associated functions. 

\smallskip \noindent 
Independently of this framework, we
have proposed in \cite{Du07,Du08}  the Taylor expansion method
for isothermal fluid problems to obtain
formally  equivalent partial differential equations at second order accuracy.
This work was  inspired by the method of equivalent equations for finite
difference schemes, in a spirit proposed  by Lerat and Peyret \cite{LP74}
and   Warming and Hyett \cite{WH74}.
In this approach, the infinitesimal variable
is simply the time step (proportional to the space step with the acoustic scaling)
and there is only one infinitesimal scale for the analysis.
This approach has been experimentaly validated with Pierre Lallemand in various contributions \cite{DL09, DL11, DL15, LD15}.
It was extended to third order accuracy for thermal and fluid problems in \cite{Du09}, 
automated at an arbitrary order in the linear case in \cite{ADGL14}, 
extended  with an external drift in \cite{DLT14} and  to relative velocities in \cite{DFG15a,DFG15b}. 
Our point of view is {\it a priori} not to consider the third and fourth order terms as approximations
of Burnett of super-Burnett equations (see {\it e.g.} \cite {AYB01, St05})   
but simply as errors of the lattice Boltzmann scheme when solving
partial differential equations of second order type. Nevertheless,
progress in the approximation of fourth order partial differential equations
have been obtained in \cite{OBD18}. 
Here, we have no particular fluid hypothesis. It can be applied to scalar or to vectorial  
\cite{Du14,Gr14}  lattice Bolzmann schemes.
This Taylor expansion method has been also revisited with the Maxwell iterations
of Yong {\it et al.} \cite{YZL16}.

\smallskip \noindent
In this contribution,
we first recall in Section 2 classic results related to the Boltzmann equation with discrete velocities. The link with
multiple relaxation times lattice Boltzmann schemes
as initially proposed by  d'Humi\`eres  \cite{DDH92} is presented in Section~3. 
Then we present the  Taylor expansion method in a very general way and the result is explicited 
with  compact formulas.
Our main result (Section~4) concerns the derivation of asymptotic partial differential equations
in the nonlinear case, up to fourth order accuracy.
The proof for  first and second orders is  presented in Section~5, with a validation for the isothermal fluid case 
in the specific case of the D2Q9 scheme. 
In Section~6, the expansion up  to third and fourth orders is detailed,
with a confrontation with  our previous work \cite{Du09}. 
The proofs of these  results need specific algebraic developments, presented in this contribution.
The linear case with constant coefficients is described in Section~7. 

\bigskip \bigskip     \noindent {\bf \large    2) \quad  Boltzmann equation with discrete velocities }    

\smallskip \noindent 
In the  space $ \,  \RR^d \, $ of dimension $\, d$, we consider a finite set
of $ \, q \, $ discrete velocities $ \, v_j  \in  {\cal V} \, $
with components $\,  v_j^\alpha \, $ for $ \, 1 \leq \alpha \leq d $. 
The unknowns of the  Boltzmann equation are the particle densities $ \, f_j $. They are 
functions of space $\, x  $, time $\, t \, $ and discrete velocities $\, v_j \, $:
%
\moneqstar 
f_j = f_j(x,\, t), \quad x \in \RR^d ,\, t \geq 0 ,\,  0 \leq j < q \, . 
\monendstar
The vector  $ \, f (x,\, t) \in \RR^q \, $ is constructed with the numbers $\, f_j(x,\, t) \, $
for $ \,  0 \leq j < q $. 
We introduce a collision vector $ \, \RR^q \ni f \longmapsto Q(f) \in  \RR^q \, $
of components $ \, Q_j(f) \, $ for $ \, 0 \leq j < q \, $ 
with given regular functions $\, Q_j $.
We introduce also a small parameter $ \, \varepsilon > 0 \, $ that can be interpreted
as a Knudsen number in the context of gas dynamics (see {\it e.g.} \cite{CC39}).

\noindent 
The Boltzmann model  with discrete velocities takes the form
\moneq \label{boltzmann-discret} 
{{\partial f_j}\over{\partial t}} + \sum_{1 \leq \alpha \leq d} v_j^\alpha \, {{\partial f_j}\over{\partial x_\alpha}} = 
{{1}\over{\varepsilon}} \, Q_j(f) ,\quad 0 \leq j < q \, . 
\monend 
It has been proposed by Carleman \cite{Ca57}, Gross \cite{Gr60}, Broadwell \cite{Br64}  and intensively developed by
Gatignol \cite{Ga75} and her co-workers.

\smallskip  \monitem {\bf  Moments }

\noindent 
Our framework concerns multi relaxation times:  
we introduce a constant invertible matrix~$ \, M \, $ 
called ``d'Humi\`eres matrix''  \cite{DDH92} in this contribution. 
This matrix defines the vector of moments $ \, m \in \RR^q \, $ by a simple product:  
\moneq \label{moments} 
 m_k \equiv \sum_{0 \leq j < q}  \, M_{k j} \,\, f_j \,.
 \monend
We introduce the number of conservation laws $ \, N \, $ ($  1 \leq N < q $) such that the  $ \, N \, $ first  moments
of the collision kernel are equal to zero:
\moneq  \label{annulation} 
\sum_{0 \leq j < q}  \, M_{k j} \,\, Q_j(f) = 0 \,, \quad \forall \, f \in \RR^q ,\,\, 0 \leq k < N \, . 
\monend
Then it is natural to divide the vector of moments  into two families: 
\moneq \label{moments-2-familles} 
m \equiv \begin{pmatrix} W \\ Y \end  {pmatrix} \, . 
\monend 
The conserved moments or macroscopic  variables $ \, W \, $ constitute 
a linear space of dimension~$ \, N \, $
and if the  Boltzmann model  with discrete velocities (\ref{boltzmann-discret}) 
is satisfied, we have  $ \, N \, $ conservation laws due to the cancellation~(\ref{annulation}):
\moneqstar
{{\partial W_k}\over{\partial t}} + \sum_{1 \leq \alpha \leq d}  \, \sum_{0 \leq j < q}  M_{k j}
\,\, v_j^\alpha \,\, {{\partial f_j}\over{\partial x_\alpha}} = 0  \,, \quad  0 \leq k < N \, .   
\monendstar
Observe that the nonconserved moments or microscopic  variables $ \, Y \, $ 
generate  a linear space of dimension $ \, q - N $.

\smallskip  \monitem {\bf  Equilibrium states }

\smallskip \noindent
We suppose that the equilibrim states $\, f^{\rm eq} \, $ defined by the conditions 
\moneqstar
Q(f^{\rm eq}) = 0 
\monendstar
are  characterized with the help of a regular nonlinear vector field $ \, \Phi \,: \RR^N \longrightarrow \RR^{q-N}  \, $
such that 
\moneq \label{f-equilibre} 
f^{\rm eq} = M^{-1} \,  \begin{pmatrix} W \\ \Phi(W) \end  {pmatrix} \, . 
\monend
In other words,  the vector field $ \, W \longmapsto Y^{\rm eq} \equiv \Phi(W) \, $ defines
the set of equilibrium states. 

\smallskip \noindent
We suppose moreover that the jacobian matrix $ \, \dd Q (f^{\rm eq}) \, $ at equilibrium
is diagonalizable with real eigenvalues and real eigenvectors. More precisely,
taking into account the hypothesis~(\ref{annulation}), 
we suppose that there exists a diagonal matrix 
\moneqstar
Z = {\rm diag} \, \Big ( {{1}\over{\tau_1}} , \, \dots ,\, {{1}\over{\tau_{q-N}}} \Big) 
\monendstar
of order $ \, q-N \, $ 
with strictly positive coefficients  $ \, \tau_\ell \, $ such that
\moneq \label{dQ-diagonalisable}  
M \, \dd Q (f^{\rm eq}) \, M^{-1}  = - \begin{pmatrix} 0 & 0 \\ 0 & Z \end {pmatrix} \,.
\monend 
In other words, the jacobian operator $ \, \dd Q (f^{\rm eq}) \, $
admits the matrix $ \, M^{-1} \, $
as a matrix of right eigenvectors and the associated eigenvalues are all non positive.

\smallskip  \monitem {\bf  Momentum-velocity operator}

\noindent 
For a linear space of dimension $ \, d $, 
we introduce the momentum-velocity operator matrix $ \, \Lambda \, $ defined by the relation 
\moneq \label{Lambda} 
 \Lambda  =  M \,\, {\rm diag} \,\Big(  \sum_{1  \leq \alpha \leq d}  v^\alpha\, \partial_\alpha \Big) \,\, M^{-1} \, . 
\monend
It is a $ \, q \times q \, $ operator matrix 
composed by first-order space differential operators.
It is obtained by conjugation of the first order advection operator $\, v . \nabla \, $
by the d'Humi\`eres matrix~$ \, M $.
The operator matrix $ \, \Lambda \, $ is nothing else than the advection operator seen in the basis of moments. 
We introduce a block decomposition  of the  momentum-velocity operator matrix associated 
to the decomposition (\ref{moments-2-familles}) of the moments.  We define 
a $ \, N \times N \, $ operator matrix $ \, A $, 
a $ \, N \times (q-N) \, $ operator matrix $ \, B $, 
a $ \,  (q-N) \times N \, $ operator matrix $ \, C \, $ and 
a $ \,  (q-N) \times (q-N) \, $ operator matrix $ \, D \, $ according to 
\moneq \label{bloc-Lambda} 
\Lambda \equiv  \begin{pmatrix} A &   B \\ C  &  D  \end  {pmatrix} \, . 
\monend
Remember that in the following, the matrices $ \, A $, $ \, B $, $ \, C \, $ and $ \, D \, $
are matrices composed with first order space  operators.

\smallskip  \monitem {\bf  Boltmann-BGK system with discrete velocities }

\noindent
We approximate now the previous discrete model with a BGK \cite{BGK54} type hypothesis:
the state~$ \, f \, $ is close to equilibrium and we approach the collision kernel
$\, Q \, $ by its first order expansion  around the  equilibrium state:
\moneqstar
Q(f) \simeq Q(f^{\rm eq}) + \dd Q(f^{\rm eq}) . (f - f^{\rm eq} ) \, .
\monendstar 
Due to the condition $ \, Q(f^{\rm eq}) = 0 $, 
we obtain with  this approximation that we call ``BGK'' in this contribution,
even it is not exactly the hypothesis done in the original article \cite{BGK54}, a new system of partial
differential equations:
\moneq \label{boltzmann-bgk-discret} 
{{\partial f_j}\over{\partial t}} + \sum_{1 \leq \alpha \leq d} v_j^\alpha \, {{\partial f_j}\over{\partial x_\alpha}} = 
{{1}\over{\varepsilon}} \, \dd Q_j(f^{\rm eq}) . (f - f^{\rm eq} ) \,, \quad 0 \leq j < q \, . 
\monend 
Observe that $ \, f^{\rm eq} \, $ is a function of the particle distribution $ \, f \, $ through the relations
 (\ref{moments}), (\ref{moments-2-familles})  and (\ref{f-equilibre}).
Taking into account the  hypothesis
(\ref{dQ-diagonalisable}), the definition  (\ref{Lambda})  and the notation (\ref{bloc-Lambda}),
we can write the  Boltmann-BGK system with discrete velocities (\ref{boltzmann-bgk-discret})
under the form of a system of two coupled  partial differential equations:
\moneq \label{systeme-bbgkg}
 {{\partial W}\over{\partial t}} + A \, W + B \, Y = 0 \,, \quad 
{{\partial Y}\over{\partial t}} + C \, W + D \, Y = -{{1}\over{\varepsilon}} \, Z \, ( Y - \Phi(W)) \, . 
\monend
%
%

\bigskip \noindent {\bf Proposition 1. \quad Fourth order Chapman-Enskog  expansion of the  Boltmann-BGK system with discrete velocities}

\noindent
As $\, \varepsilon \, $ tends to zero,
the conserved moments $ \, W \, $ of the Boltmann-BGK system with discrete velocities~(\ref{systeme-bbgkg})
satisfy formally the following asymptotic  system of partial differential equations:
\moneq \label{edp-W-temps-continu} 
 {{\partial W}\over{\partial t}} + \Gamma_1(W) + \varepsilon \, \Gamma_2(W)
+ \varepsilon^2 \, \Gamma_3(W) +  \varepsilon^3 \, \Gamma_4(W) = {\rm O}(\varepsilon^4) \, .
\monend
The nonlinear operators are related to the expansion of the microscopic moments:
\moneq \label{dvpt-Y-temps-continu} 
Y = \Phi(W) + \varepsilon \, Z^{-1} \, \Psi_1(W) +  \varepsilon^2 \, Z^{-1} \,  \Psi_2(W) +  \varepsilon^3 \, Z^{-1} \, \Psi_2(W) + {\rm O}(\varepsilon^4)
\monend
and we have the interlaced relations
\moneq \label{formules-temps-continu} \left \{ \begin {array}{rl}
\Gamma_1 (W)  & \!\!\! =     A \, W + B \, \Phi(W) \\
\Psi_1 (W) & \!\!\! =   \dd \Phi(W) .  \Gamma_1 (W)  - \big( C \, W + D \, \Phi(W)  \big) \\ 
\Gamma_2 (W) & \!\!\! =  B \, Z^{-1}  \, \Psi_1 (W) \\
\Psi_2 (W) & \!\!\! =   Z^{-1}  \,  \dd \Psi_1 (W) .   \Gamma_1 (W) + \dd \Phi(W) .  \Gamma_2 (W) 
- D \, Z^{-1}  \,  \Psi_1 (W)   \\ 
\Gamma_3(W)   & \!\!\! = B \, Z^{-1}   \, \Psi_2 (W)  \\
\Psi_3 (W) & \!\!\! = Z^{-1}  \, \dd \Psi_1 (W) .  \Gamma_2 (W)  +  \dd \Phi(W) .  \Gamma_3(W) -  D \, Z^{-1}  \, \Psi_2 (W)  \\
&   + Z^{-1}  \, \dd \Psi_2 (W) .  \Gamma_1 (W) \\
\Gamma_4(W)   & \!\!\! =      B \, Z^{-1}  \, \Psi_3 (W)\, . 
\end {array} \right. \monend

\smallskip \noindent
In applications to fluid dynamics, the first order term $ \, \Gamma_1 (W) \, $ is associated to
the Euler equations of gas dynamics and  the second order term $ \, \Gamma_2 (W) \, $ to the viscous terms
of the Navier-Stokes equations.
The  operators $ \, \Gamma_3 (W) \, $ and $ \, \Gamma_4 (W) \, $ represent 
the Burnett and  super-Burnett terms  respectively. 
With this Proposition 1, we essentially reformulate the classic Chapman Enskog expansion  \cite{CC39}
in the discrete velocity case, as done {\it e.g.} in the book of Gatignol \cite{Ga75}.

\monitem Proof of Proposition 1. 

\noindent
First  we inject the representation (\ref{dvpt-Y-temps-continu}) inside the first equation of (\ref{systeme-bbgkg}).
It comes
\moneqstar
{{\partial W}\over{\partial t}} + A \, W + B \, \big( \Phi(W)
+ \varepsilon \, Z^{-1} \, \Psi_1(W) +  \varepsilon^2 \, Z^{-1} \,  \Psi_2(W) +  \varepsilon^3 \, Z^{-1} \, \Psi_2(W) \big) 
= {\rm O}(\varepsilon^4) \,  
\monendstar
then after rearranging the terms, 
\moneqstar 
{{\partial W}\over{\partial t}} + A \, W  + B \, \Phi(W) +  \varepsilon \,  B \,   Z^{-1} \, \Psi_1(W)
+  \varepsilon^2 \,  B \,   Z^{-1} \, \Psi_2(W) +  \varepsilon^3 \,  B \,   Z^{-1} \, \Psi_3(W) = {\rm O}(\varepsilon^4) \,. 
\monendstar
The confrontation of this partial differential equation with the Ansatz (\ref{edp-W-temps-continu}) shows that 
%
%
\moneqstar  \left \{ \begin {array}{rl}
\Gamma_1 (W)  & \!\!\! =     A \, W + B \, \Phi(W) \\
\Gamma_2 (W) & \!\!\! =  B \, Z^{-1}  \, \Psi_1 (W) \\
\Gamma_3(W)   & \!\!\! = B \, Z^{-1}   \, \Psi_2 (W)  \\
\Gamma_4(W)   & \!\!\! =      B \, Z^{-1}  \, \Psi_3 (W) 
\end {array} \right. \monendstar
as proposed in (\ref{formules-temps-continu}). 

\noindent
Secondly, we write the second equation of (\ref{systeme-bbgkg}) under the form
\moneqstar
Y = \Phi - \varepsilon \,  Z^{-1}  \, ( \partial_t Y + C \, W + D \, Y ) \, . 
\monendstar 
We observe that

\smallskip \noindent $ -\partial_t Y = \dd Y\, . \,(- \partial_t W)  $

\noindent \qquad  $ \,\,\,\,
= \dd \big( \Phi(W) + \varepsilon \, Z^{-1} \, \Psi_1(W) +  \varepsilon^2 \,  Z^{-1} \, \Psi_2(W) \big)
\,.\, \big( \Gamma_1 +  \varepsilon \, \Gamma_2(W) + \varepsilon^2 \, \Gamma_3(W) \big) + {\rm O}(\varepsilon^3)  $ 

\noindent \qquad  $ \,\,\,\,
= \dd \Phi(W) \,.\,\Gamma_1 +  \varepsilon \, \big( Z^{-1} \, \dd \Psi_1(W)\,.\,\Gamma_1 +  \dd \Phi(W) \,.\,\Gamma_2 \big)  $

\qquad \qquad $ 
+  \varepsilon^2 \, \big( Z^{-1} \, \dd \Psi_2(W)\,.\,\Gamma_1 + Z^{-1} \, \dd \Psi_1(W)\,.\,\Gamma_2  +  \dd \Phi(W) \,.\,\Gamma_3 \big)
+  {\rm O}(\varepsilon^3) \, . $

\smallskip \noindent
We report this expression inside the value of $ \, Y \, $:

\smallskip \noindent
$ Y = \Phi +  \varepsilon \,  Z^{-1} \, (  -\partial_t Y -  C \, W - D \, Y ) $

\noindent \quad $ \, =
\Phi +  \varepsilon \,  Z^{-1} \, \big[  \dd \Phi(W) \,.\,\Gamma_1 +  \varepsilon \, ( Z^{-1} \, \dd \Psi_1(W)\,.\,\Gamma_1 +  \dd \Phi(W) \,.\,\Gamma_2 )  $

\noindent \qquad \qquad \qquad \quad    $
+  \varepsilon^2 \, ( Z^{-1} \, \dd \Psi_2(W)\,.\,\Gamma_1 + Z^{-1} \, \dd \Psi_1(W)\,.\,\Gamma_2  +  \dd \Phi(W) \,.\,\Gamma_3 ) $

\noindent \qquad \qquad \qquad \quad   $
- \, C \, W  - D \, \big( \Phi(W) + \varepsilon \, Z^{-1} \, \Psi_1(W) +  \varepsilon^2 \,   Z^{-1} \, \Psi_2(W) \big) \big]  + {\rm O}(\varepsilon^4)  $ 

\noindent \quad $ \, =
\Phi +  \varepsilon \,  Z^{-1} \, ( \dd \Phi(W) \,.\,\Gamma_1 -C \, W - D \, \Phi )
 +  \varepsilon^2 \,  Z^{-1} \, (  Z^{-1} \, \dd \Psi_1(W)\,.\,\Gamma_1 +  \dd \Phi(W) \,.\,\Gamma_2 - D \,  Z^{-1} \, \Psi_1 ) $

\noindent \qquad \quad   $
+  \varepsilon^3 \,  Z^{-1} \, (   Z^{-1} \, \dd \Psi_2(W)\,.\,\Gamma_1 + Z^{-1} \, \dd \Psi_1(W)\,.\,\Gamma_2  +  \dd \Phi(W) \,.\,\Gamma_3 
- D \,  Z^{-1} \, \Psi_2 )  + {\rm O}(\varepsilon^4) \, . $

\smallskip \noindent
Comparatively to  the relation (\ref{dvpt-Y-temps-continu}), we obtain
\moneqstar  \left \{ \begin {array}{l}
 \Psi_1 (W) = \dd \Phi(W) \,.\,  \Gamma_1  - ( C \, W + D \, \Phi ) \\
 \Psi_2 (W) = Z^{-1}  \,  \dd \Psi_1 (W) \,.\, \Gamma_1 + \dd \Phi(W) \,.\,  \Gamma_2 (W) - D \, Z^{-1}  \,  \Psi_1 (W)  \\
 \Psi_3 (W) = Z^{-1}  \, \dd \Psi_1 (W) \,.\,  \Gamma_2  +  \dd \Phi(W) \,.\,  \Gamma_3 -  D \, Z^{-1}  \, \Psi_2 + Z^{-1}  \, \dd \Psi_2 (W) \,.\,  \Gamma_1   
\end {array} \right. \monendstar
and the last relations of  (\ref{formules-temps-continu}) are explicited. The proof is completed. \hfill $\square$

\newpage 
\bigskip \bigskip     \noindent {\bf \large    3) \quad  Multiple Relaxation Times Lattice Boltzmann schemes }    

\smallskip \noindent
We introduce now   
a regular cartesian lattice $ \, {\cal L} \, $ composed by vertices $ \, x  \, $
separated by distances that are simple expressions of the space step $ \, \Delta x $. 
A discrete time~$ \, t \, $ is supposed to be an integer  multiple of a time step $ \, \Delta t > 0 $.
The unknowns of a lattice Boltzmann
scheme  with $ \, q \, $ discrete velocities are the particle densities $ \, f_j $. 
They are now functions of discrete space $\, x $, discrete time $\, t \, $ and discrete velocities $\, v_j \, $
for $ \, 0 \leq j < q $: 
\moneqstar 
f_j = f_j(x,\, t), \quad x \in    {\cal L} , \, t = n \, \Delta t ,\, n \in \NN ,\,  v_j \in  {\cal V} \, . 
\monendstar

\smallskip \noindent 
The discrete velocity set $ \,  {\cal V} \, $ does not depend on space or time.

\smallskip \noindent
Multiple relaxation times lattice Boltzmann scheme considered in this contribution consists
in a simple Lie-Trotter splitting applied to the Boltzmann-BGK model with discrete velocities~(\ref{boltzmann-bgk-discret}).
First a relaxation step $ \, f(t) \longrightarrow R(\Delta t)  \, f(t) \equiv f^* \, $ for the approximation
of the system of ordinary differential equations
\moneqstar
{{\partial f}\over{\partial t}}  = {{1}\over{\varepsilon}} \, \dd Q(f^{\rm eq}) . (f - f^{\rm eq} )
\monendstar
and secondly a transport  step $ \,  f^*   \longrightarrow T(\Delta t) \, f^* \equiv f(t+\Delta t) \, $
for the free advection  
\moneqstar 
{{\partial f_j}\over{\partial t}} + \sum_{1 \leq \alpha \leq d} v_j^\alpha \, {{\partial f_j}\over{\partial x_\alpha}} = 0 \,, \quad 0 \leq j < q \, .
\monendstar 
%

\smallskip  \monitem {\bf  Relaxation step }

\noindent
Due to the formulation (\ref{systeme-bbgkg}), 
the relaxation scheme is simply described in the moment representation and we have
\moneqstar
 {{\partial W}\over{\partial t}}  = 0 \,, \quad 
{{\partial Y}\over{\partial t}}   =  -{{1}\over{\varepsilon}} \, Z \, ( Y - \Phi(W)) \, . 
 \monendstar 
%
%
After relaxation, the vector of moments $ \, m \, $ is transformed into a new vector
$ \, m^* \, $: 
\moneq \label{m-star} 
m^* \equiv  \begin{pmatrix} W^* \\ Y^* \end  {pmatrix} \, . 
\monend 
The conserved moments are not modified during this relaxation step: 
\moneq \label{relaxation-W} 
 W^* = W \, . 
\monend
As remarked in \cite{Du08}, we use an explicit first order Euler scheme for the microscopic moments:
\moneqstar
{{Y^* - Y}\over{\Delta t}} =   -{{1}\over{\varepsilon}} \, Z \, ( Y - \Phi(W)) \,, 
\monendstar 
{\it id est }
\moneq \label{relaxation-Y} 
 Y^* = Y + S \, ( \Phi(W) - Y ) \,  
\monend 
with 
\moneq \label{def-matrice-S}  
S = {{\Delta t}\over{\varepsilon}} \, Z \, .
\monend
Usually, as pointed by Lallemand and Luo  \cite{LL00}, the relaxation matrix $ \, S \, $  is a diagonal matrix, in coherence
with the hypothesis concerning the diagonalization of the operator $\, \dd Q(f^{\rm eq}) \, $
presented in the previous  section:
\moneqstar
S = {\rm diag} \, \big ( s_1 , \, \dots ,\, s_{q-N}  \big) \, . 
\monendstar
Moreover, the  matrix $ \, S \, $  is a fixed invertible square matrix of dimension $ \, (q-N)\times(q-N) $. 
In other words, the time step $ \, \Delta t \, $ is if the order of the Knudsen number
$ \, \varepsilon \, $ in our asymptotic study. 
 Finally the particle distribution $ \, f^* \, $ after relaxation is given according to 
\moneq \label{f-star} 
f^* = M^{-1} \, m^*  \, . 
\monend 
%

\smallskip  \monitem {\bf  Necessary conditions for stability }

\noindent
We deduce from the relation (\ref{relaxation-Y})
\moneqstar  
 Y^* - \Phi(W) = ({\rm I} -S) \, (Y - \Phi(W) ) \,  
\monendstar 
and a natural stability condition takes the form $ \, | | \, {\rm I} -S \, | | \leq 1 $, {\it id est}
\moneqstar
0 \leq s_j \leq 2 \,, \quad 0 \leq j < q-N \, . 
\monendstar 
The case of over-relaxation $ \, 1 < s_j < 2 \, $ is widely used
when implementing  lattice Boltzmann schemes to take into account very law viscosities in fluid mechanics.
We will emphasize  this point in the next section.

\smallskip  \monitem {\bf  Advection step }

\noindent
In the second step of advection, the 
method of characteristics is applied in the very particular case when it is  exact.
In one time step each particle distribution $ \, f_j \, $ is exactly transported
from the note $ \, x - v_j \, \Delta t \, $ of the lattice $ \, {\cal L} \, $ to the vertex $ \, x \in  {\cal L} $.
In other words, the advection step  occurs with a Courant-Friedrichs-Lewy number $ \, cfl \equiv 1 \, $
for all discrete velocities $ \, v_j $. 

\smallskip  \monitem {\bf  Lattice Boltzmann scheme as a Lie-Trotter splitting scheme }

\noindent
The iteration of the lattice Boltzmann scheme is finally given by the relation
\moneq \label{lie-trotter} 
f(t+\Delta t) = T(\Delta t) \, R(\Delta t) \, f(t) \, 
\monend
and it is usefull to explicit it for each component: 
\moneq \label{iteration-f} 
 f_j(x,\, t+\Delta t) \,=\, f_j^*(x - v_j \, \Delta t , \, t)  \,, \quad v_j \in {\cal V}
\,,  \,\, x \in  \cal{L} \, . 
\monend
Observe that this framework is very general.
The so-called ``BGK version'' of the lattice Boltzmann scheme  \cite{BC98, HSB91, QHL92} 
corresponds to the choice $ \, M = {\rm I} \, $ and $ \, Z = {{1}\over{\tau}} \,  {\rm I} $.
Usual multi-relaxation  lattice Boltzmann schemes  \cite{DGL16, DL09,DDH92, LL00}
suppose in general some orthogonality property for the d'Humi\`eres matrix. 
This paradigm  permits also  the introduction of two particle distributions
\cite{ACS93, DGL16, GZS07, KSC13}              
or eventually more with the vectorial schemes developed by Graille  \cite{Gr14}
or a variant proposed in \cite{Du14}.         
In those last cases, the numbering $ \, j \longmapsto v_j \, $ of the velocities
is simply non injective. 

\smallskip \noindent
It is naturally possible to think to a second order accurate splitting scheme. A natural idea is the Strang splitting \cite{St68}:
\moneqstar
f(t+\Delta t) =  R \Big( {{\Delta t}\over{2}} \Big)  \,\,  T(\Delta t)  \,\, R \Big( {{\Delta t}\over{2}} \Big) \, f(t) \, .
\monendstar 
In this case, as put in evidence by Dellar \cite{De13}, one can write $ \, \widetilde{f}(t) \equiv  R \big( {{\Delta t}\over{2}} \big)^{-1} \, f(t) \, $
and the  previous Strang scheme can be written as 
\moneqstar
 \widetilde{f}(t+\Delta t) =  T(\Delta t)  \,\, R \Big( {{\Delta t}\over{2}} \Big)  \, R \Big( {{\Delta t}\over{2}} \Big) \,  \widetilde{f}(t)  \, .
\monendstar 
This scheme, up to the shift of one half relaxation step, is identical to the Lie-Trotter splitting~(\ref{lie-trotter}), except that the
product $ \,  R \big( {{\Delta t}\over{2}} \big)  \, R \big( {{\Delta t}\over{2}} \big) \, $ plays the role of $ \, R(\Delta t) $.
A first possibility is to divide the time step $ \, \Delta t \, $ by a factor $ \,2 \, $ in the relations
(\ref{relaxation-Y})(\ref{def-matrice-S}):
\moneqstar
\overline {Y} = Y + {1\over2} \, S \, (\Phi - Y) \,,  \quad Y^* = \overline {Y} + {1\over2} \, S \, (\Phi -  \overline Y) \, . 
\monendstar 
Then  
\moneqstar
Y^* = \Big( 1 - S + {1\over4} \, S^2 \Big) \, Y + \Big( S - {1\over4} \, S^2 \Big) \, \Phi(W) 
\monendstar 
and over-relaxation is no more possible because $ \, s_j -  {1\over4} \, s_j^2 \leq 1 \, $ for all $\, s_j \in \RR $. 
So such a numerical scheme is not used by the engineers at our knowledge.

\smallskip \noindent
A second possibility is to enforce the relation $ \,  R \big( {{\Delta t}\over{2}} \big) = \sqrt{ R(\Delta t) } $, {\it id est}
\moneqstar
\overline {Y} = \sqrt{ {\rm I} - S} \,\, Y + \big(  {\rm I} -  \sqrt{ {\rm I} - S} \big)  \, \Phi(W)  \, . 
\monendstar 
But this relation has a sense only if $ \, s_j \leq 1 \, $ and over-relaxation is again excluded. 
  
\smallskip \noindent
An other idea is to exchange the roles of advection and relaxation in the Strang scheme:
\moneqstar
f(t+\Delta t) =  T \Big( {{\Delta t}\over{2}} \Big)  \,\,  R(\Delta t)  \,\, T \Big( {{\Delta t}\over{2}} \Big) \, f(t) \, 
\monendstar 
or to use more elaborated splittings as proposed by Drui {\it et al.} \cite {DFHN19}.
In this case, the transport step is nomore at a Courant number equal to unity and an interpolation
procedure, costly in terms of numerical viscosity or complex to implement as the semi-Lagrangian method \cite{Fl19},
is not usually chosen by experts developping lattice Boltzmann schemes.

\smallskip \noindent
In conclusion of this sub-section, the naive Lie-Trotter splitting scheme (\ref{lie-trotter})
gives very good numerical results. It is important to put in evidence precise qualities and defects of this scheme
through the following  formal asymptotic analysis.


\bigskip \noindent {\bf Proposition 2. \quad Exponential form of  lattice Boltzmann schemes}  

\noindent 
Consider a lattice Boltzmann scheme as defined at the previous section
by the relations
(\ref{moments}),    (\ref{moments-2-familles}),    (\ref{f-equilibre}), (\ref{Lambda}),  (\ref{bloc-Lambda}), 
  (\ref{m-star}),  (\ref{relaxation-W}),   (\ref{relaxation-Y}),  (\ref{f-star}) and (\ref{iteration-f}). 
Then we have an exponential form of the discrete iteration (\ref{iteration-f}): 
\moneq \label{exponentielle-Lambda} 
m  (x, t + \Delta t) =  {\rm exp} ( - \Delta t \, \Lambda ) \,\,  m^*(x ,\, t) \, .    
\monend
%

\monitem Proof of Proposition 2. 

\noindent We have simply a long sequence of identities: 

 \smallskip  \noindent $ \displaystyle 
\, m_k (x, t + \Delta t) = \sum_{j} M_{k j} \, f_j(x, t + \Delta t) $ \hfill due to (\ref{moments})

\qquad \qquad \qquad  $ \, \displaystyle  = 
\sum_{j} M_{k j} \, f_j^* (x-v_j \, \Delta t ,\, t) $  \hfill due to (\ref{iteration-f}) 

\qquad \qquad \qquad  $ \, \displaystyle  = 
\sum_{j \, \ell} M_{k j} \, (M^{-1})_{_{\scriptstyle \! j \ell}} \, m_\ell^*(x-v_j \, \Delta t ,\, t) $ \hfill due to (\ref{f-star})   

\qquad \qquad \qquad  $ \, \displaystyle  = 
\sum_{j \, \ell} M_{k j} \, (M^{-1})_{_{\scriptstyle \! j \ell}} \, \sum_{n=0}^{\infty} 
{{1}\over{n !}} \big(-\Delta t \sum_{\alpha}  v_j^\alpha \, \partial_\alpha \big)^n  \,  m_\ell^*(x ,\, t) $ 
\hfill Taylor expansion    

\qquad \qquad \qquad  $ \, \displaystyle  = 
\sum_{\ell} \sum_{n=0}^{\infty} 
{{1}\over{n !}} \, \sum_{j} M_{k j} \,  \Big( \! -\Delta t \, \sum_{\alpha}  v_j^\alpha \, \partial_\alpha \Big)^n \, (M^{-1})_{_{\scriptstyle \! j \ell}} 
\, m_\ell^*(x ,\, t) $
\hfill elementary algebra

\qquad \qquad \qquad  $ \, \displaystyle  = 
\sum_{\ell} \sum_{n=0}^{\infty} 
{{1}\over{n !}} \,   (-\Delta t)^n \,  \bigg( \Big( M  \, \big( \sum_{\alpha}  v^\alpha  \, \partial_\alpha  \big) \, M^{-1}  \Big)^{\!\! n} \bigg)_{\scriptstyle \!\! k \ell}
\,\,  m_\ell^*(x ,\, t) $

\hfill because $ \, M \, $ is a constant matrix

\qquad \qquad \qquad  $ \, \displaystyle  = 
\sum_{\ell} \sum_{n=0}^{\infty} 
{{1}\over{n !}} \,   (-\Delta t)^n \,  \big( \Lambda^n \big)_{\scriptstyle \! k \ell}
\,\,  m_\ell^*(x ,\, t) $
\hfill due to (\ref{Lambda})

\qquad \qquad \qquad  $ \, \displaystyle  = 
\sum_{\ell}  \bigg[  \, \sum_{n=0}^{\infty} 
{{1}\over{n !}}  \,  \big(- \Delta t \, \Lambda \big)^n_{k \ell}  \, \bigg]  \,\,  m_\ell^*(x ,\, t) $
\hfill elementary algebra     

\qquad \qquad \qquad  $ \, \displaystyle  = 
\sum_{\ell} \,   {\rm exp} ( - \Delta t \, \Lambda )_{k \ell}   \, \,  m_\ell^*(x ,\, t) $ 
\hfill elementary algebra     

\qquad \qquad \qquad  $ \, \displaystyle  = 
\Big(  {\rm exp} ( - \Delta t \, \Lambda ) \,\, m^*(x ,\, t) \Big)_{\! \! k} $
\hfill elementary algebra     

and the proof is completed. \hfill $\square $

\smallskip 
We can see the relation (\ref{exponentielle-Lambda}) as a discrete form of the Duhamel formula
for ordinary differential equations.
The idea of introducing a formal exponential expansion is also present in
the work of Boghosian and Coveney in the BGK case \cite{BC98}. 

\smallskip  \monitem {\bf  Momentum-velocity operator for a D2Q9 lattice Boltzmann scheme } 

\begin{figure}    [H]  \centering 
\vspace{-1.5 cm} 
\centerline { \includegraphics[width=.45  \textwidth, angle=0] {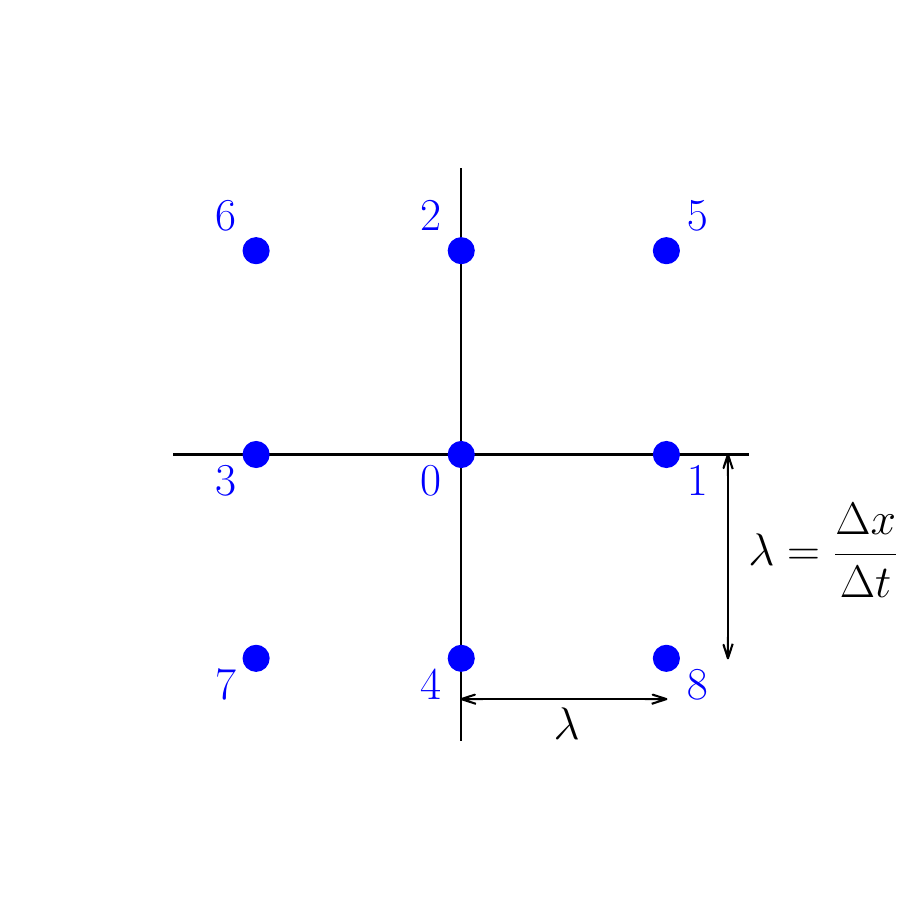}}
\vspace{-1.5 cm} 
\caption{D2Q9 lattice Boltzmann scheme } 
\label{fig-d2q9} \end{figure}

\noindent 
We can concretize the  matrix $ \, \Lambda \, $ introduced in (\ref{Lambda})
with the popular two-dimensional ``D2Q9'' scheme.
A set of nine velocities is defined by the Figure~\ref{fig-d2q9}. 
%
%
To fix the ideas, we use the moment matrix presented in Lallemand and Luo \cite{LL00}: 
\moneq \label{M-d2q9}
M_{D2Q9}  =    \left[  \!\!   \begin{array} {ccccccccc}
 1  \!&\!  1  \!&\!  1  \!&\!  1  \!&\!  1  \!&\!  1  \!&\!  1  \!&\!  1  \!&\!  1 \\ 
   0 \!&\! \lambda  \!&\!  0  \!&\!  -\lambda  \!&\!  0  \!&\!   \lambda  \!&\!   -\lambda   \!&\!  -\lambda  \!&\!  \lambda  \\ 
 0  \!&\!  0   \!&\! \lambda  \!&\!  0    \!&\!   -\lambda   \!&\!   \lambda   \!&\!   \lambda  \!&\!  -\lambda  \!&\!  -\lambda \\ 
-4 \lambda^2    \!&\!  -\lambda^2   \!&\!   -\lambda^2   \!&\!   -\lambda^2   \!&\!   -\lambda^2  \!&\!   2 \lambda^2   \!&\!  2 \lambda^2   \!&\!  2 \lambda^2   \!&\!  2 \lambda^2 \\ 
  0    \!&\!   \lambda^2    \!&\!  -\lambda^2    \!&\!   \lambda^2   \!&\!   -\lambda^2  \!&\!  0  \!&\!  0  \!&\!  0  \!&\!  0 \\ 
 0  \!&\!  0  \!&\!  0  \!&\!  0  \!&\!  0   \!&\!   \lambda^2   \!&\!   -\lambda^2    \!&\!   \lambda^2   \!&\!   -\lambda^2 \\ 
 0   \!&\! -2 \lambda^3  \!&\!  0   \!&\!  2 \lambda^3  \!&\!  0   \!&\!  \lambda^3   \!&\!   -\lambda^3   \!&\!   -\lambda^3   \!&\!    \lambda^3 \\
 0  \!&\!  0  \!&\! -2 \lambda^3  \!&\!  0  \!&\!  2 \lambda^3    \!&\!   \lambda^3    \!&\!   \lambda^3    \!&\!  -\lambda^3    \!&\!  -\lambda^3 \\ 
4 \lambda^4  \!&\! -2 \lambda^4  \!&\! -2 \lambda^4  \!&\! -2 \lambda^4  \!&\! -2 \lambda^4  \!&\!   \lambda^4   \!&\!  \lambda^4   \!&\!  \lambda^4  \!&\!   \lambda^4 
\end{array} \!\! \right] \, . 
\monend  
Observe that other choices are possible with the Hermite moments developed
in \cite{PH06,SHC06}. 
The moments are associated with the lines of the matrix (\ref{M-d2q9}).
They are denominated (in this order) by
$ \, \rho ,\, J_x ,\, J_y ,\, \varepsilon, \, x_x ,\, x_y ,\, q_x ,\, q_y ,\, h $. 
The lines of this invertible matrix are chosen orthogonal.
For isothermal flows, the conserved moments
\moneq \label{d2q9-W} 
W =  \big( \rho \,,\,  J_x  \,,\,  J_y \big)^{\rm t}
\monend 
 correspond to the three first lines of the d'Humi\`eres matrix
 (\ref{M-d2q9}).
The non-conserved moments complete the family: 
\moneq  \label{d2q9-Y}  
Y =  \big( \varepsilon \,,\, x_x \,,\, x_y \,,\, q_x \,,\, q_y  \,,\, h \big)^{\rm  t} \, . 
\monend 
They are  linked to the six last lines of the matrix defined in (\ref{M-d2q9}).
With the velocities defined in Figure~\ref{fig-d2q9} and the moment matrix  (\ref{M-d2q9}),
we construct without difficulty the
operator matrix  $ \, \Lambda \, $ for the thermal  D2Q9 scheme:
\moneq \label{Lambda-d2q9}
\Lambda_{D2Q9}  =    \left[ \,\, \begin{array}{|ccc|cccccc|} \hline 
 0 & \!\! \partial_x & \!\! \partial_y &   0  &   0 & \!\!   0 & \!\! 0  &    0  &    0  \\ 
{{2\lambda^2}\over{3}} \, \partial_x & \!\!  0 & \!\!  0 &   {1\over6} \, \partial_x & 
  {1\over2} \, \partial_x  &     \!\!  \partial_y  & \!\!   0  &   0  &   0 \\
{{2\lambda^2}\over{3}} \, \partial_y & \!\!  0 & \!\!  0 &   {1\over6} \, \partial_y & 
   -{1\over2} \, \partial_y  &    \!\!  \partial_x  & \!\!   0  & 0 &   0
\\  \vspace{-.4 cm} & & & & & & & & \\ \hline  \vspace {-.4 cm} & & & & & & & &  \\
 0 &  \!\!  \lambda^2 \, \partial_x & \!\!   \lambda^2 \, \partial_y & 
   0 &    0 & \!\!    0 & \!\!    \partial_x &     \partial_y &    0  \\ 
  0 &  \!\!   {{\lambda^2}\over{3}} \, \partial_x &  \!\!   -{{\lambda^2}\over{3}} \, \partial_y &   0 &   0 & \!\! 
  0 & \!\!   - {1\over3} \, \partial_x &    {1\over3} \, \partial_y   &   0  \\ 
 0 &  \!\!   {{2}\over{3}}  \lambda^2 \,\partial_y &  \!\!    {{2}\over{3}}   \lambda^2 \, \partial_x & \!\!     0 &    0 & \!\! 
   0 & \!\!     {1\over3} \, \partial_y &  \!\!    {1\over3} \, \partial_x &    0   \\ 
  0 & \!\!    0 & \!\!    0 &      {{\lambda^2}\over{3}} \, \partial_x &   - \lambda^2 \, \partial_x & 
\!\!    \lambda^2 \, \partial_y & \!\!    0 &    0 &    {1\over3} \, \partial_x   \\ 
  0 & \!\!    0 & \!\!    0 &      {{\lambda^2}\over{3}} \, \partial_y &     \lambda^2 \, \partial_y & 
   \lambda^2 \, \partial_x & \!\!   0 &      0 & \!\!     {1\over3} \, \partial_y  \\ 
0 & \!\!    0 & \!\!   0 &   0 &   0 &   0 &   \lambda^2 \, \partial_x &   \lambda^2 \, \partial_y & 0   \\
\hline    \end{array}  \,\,   \right] \, . 
\monend

\smallskip  \noindent 
In the relation (\ref{Lambda-d2q9}), we have put in evidence the  block decomposition (\ref{bloc-Lambda}) 
for the D2Q9 scheme with the conserved moments precised in (\ref{d2q9-W}).
Even if the space differential operators commute, the matrices that compose the  momentum-velocity operator matrix
does not commute! We have the very important but elementary structure  
\moneqstar 
 \Lambda^2 \equiv  \begin{pmatrix} A_2 & B_2 \\ C_2 & D_2  \end  {pmatrix} = 
 \begin{pmatrix} A &   B \\ C  &  D  \end  {pmatrix}  \,\,  \begin{pmatrix} A &   B \\ C  &  D  \end  {pmatrix}  =
 \begin{pmatrix} A^2 + B \, C  &  A \, B + B \, D \\ C \, A + D \, C   & C \, B + D^2  \end  {pmatrix} \,  
\monendstar 
and to fix the ideas  
\moneqstar 
 \Lambda^3 \equiv  \begin{pmatrix} A_3 & B_3 \\ C_3 & D_3  \end  {pmatrix} = 
 \begin{pmatrix} A_2 &   B_2 \\ C_2  &  D_2  \end  {pmatrix}  \,\,  \begin{pmatrix} A &   B \\ C  &  D  \end  {pmatrix}  =
 \begin{pmatrix} A_2 \, A + B_2 \, C  &  A_2 \, B + B_2 \, D \\ C_2 \, A + D_2 \, C   & C_2 \, B + D_2 \, D \end  {pmatrix} \,  . 
\monendstar 

\bigskip \bigskip   \noindent {\bf \large   4) \quad  Taylor expansion method}    

\smallskip 
In this section, we revisit the Taylor expansion method:
we specify  the hypotheses, express our general result and begin the proof with the order zero.

\smallskip \monitem {\bf  Hypotheses for a formal expansion }

\noindent
The formal Taylor expansion supposes some precise hypotheses.
First, we adopt the acoustic scaling: the ratio $ \, \lambda \equiv {{\Delta x}\over{\Delta t}} \, $
is supposed to be constant in all this work.
%
Moreover, the ratio $ \, {{\Delta t}\over{\varepsilon}} \, $  remains constant in our asymptotic analysis.
Then the relaxation matris $ \, S \, $ is fixed and invertible;
it is also the case for the matrix $ \, S^{-1} \, $ proportional to the matrix $ \,  Z^{-1} \, $ who
plays an important role in the relations (\ref{formules-temps-continu}) of the Chapman-Enskog expansion
(\ref{edp-W-temps-continu})(\ref{dvpt-Y-temps-continu})(\ref{formules-temps-continu}) 
for the  continuous in time lattice Boltzmann-BGK system.

\smallskip \noindent
A  unintuitive result has been discovered by  H\'enon \cite{He87}.
Due to the fact that the  fully discretized in time  lattice Boltzmannn scheme (\ref{iteration-f})
is substantially different from  the   continuous in time Boltzmann-BGK system with discrete velocities
(\ref{boltzmann-bgk-discret}), the asymptotic expansions have similarities and differences.
In particular, the lattice Boltzmannn scheme 
put in evidence what we call the  H\'enon  matrix $ \, \Sigma \, $ in this contribution.
Roughtly speaking, it plays a role analogous to the matrix $ \,  Z^{-1} \, $ in the relations (\ref{formules-temps-continu}).
But it is defined by 
\moneq \label{Henon} 
 \Sigma \equiv  S^{-1} - {1\over2} \, {\rm I} \, .
\monend
This matrix emerges from the very classic second order analysis presented thereafter at  Proposition 5
and detailed in the relations (\ref{edp-ordre-2})(\ref{Gamma-2}). For applications to fluid dynamics,
this matrix is closely related to viscosities, as explained in Section 5 for the D2Q9 example 
in relations (\ref{Sigma-d2q9}) and (\ref{viscosites-d2q9}).
Then some values of the H\'enon matrix $ \, \Sigma \, $ are chosen as small as possible in order
to simulate flows with high Reynolds number.
In consequence, over-relaxation is a mandatory practice for lattice Boltzmann schemes
applied to  high Reynolds number flows. 
Recall that this matrix remains fixed in our analysis.


\smallskip \noindent
Secondly, we suppose also  that we can differentiate all the expansions relative to space and/or time.
This hypothesis is mathematically absolutly non trivial and expresses
that all the formal expansions should take place in very regular functional spaces. 
Thirdly, we use the notation
$ \, \dd \zeta(W) . \xi \, $ 
for the action of the differential of some regular function
$ \, W \longmapsto  \zeta(W) \, $ againts a test vector $ \, \xi $.
In particular, 
\moneq \label{notation-gamma-j} 
\gamma_j (W) \equiv   \dd \Phi(W) .  \Gamma_j (W) \, , \quad 1 \leq j \leq 3 \, . 
\monend 
For second order derivatives, we introduce
\moneq \label{notation-d2-psi-gamma1}
\partial^2 \Psi .  \Gamma_1  \equiv  \dd \, ( \dd  \Psi (W) . \, \Gamma_1 ) . \, \Gamma_1  =   \dd^2 \Psi (W) . (\Gamma_1 ,\, \Gamma_1) 
+  \dd \Psi (W) . \dd  \Gamma_1 (W)  .  \, \Gamma_1 \, . 
\monend   
The asymptotic analysis occurs for a time step $ \,  \Delta t \, $ (or a space step $ \, \Delta x$) tending to zero. 
Emerging partial differential equations for the conserved variables are denoted as 
\moneq   \label{edp-ordre-4} 
\partial_t W + \Gamma_1 (W) + \Delta t \, \Gamma_2 (W) + \Delta t^2 \, \Gamma_3 (W)  + \Delta t^3 \, \Gamma_4 (W) = {\rm O}(\Delta t^4) \, . 
\monend  
The vector $\, \Gamma_j (W) \, $ belongs in $ \, \RR^N $. It is obtained after $j$ space derivations of the conserved moments $ \, W \, $ and 
of the equilibrium vector $ \, \Phi(W) $.
For this reason, we will speak in the following of {\it e.g. } $ \, \Gamma_1 (W) \, $ as ``first order term'', even if it is a term
of order zero relative to the infinitesimal $ \, \Delta t $. 
For non-conserved moments (or microscopic variables), we suppose   
\moneq \label{variables-microscopiques} \left \{ \begin {array}{rl}
Y & \!\!\! = \Phi(W) +  \,\,  S^{-1} \, \big( \, \Delta t \,\, \Psi_1 (W)  + \Delta t^2 \,\, \Psi_2 (W) + \Delta t^3 \,\, \,\Psi_3 (W) \, \big) 
+ {\rm O}(\Delta t^4) \\ 
Y^* & \!\!\! = \Phi(W) +  \,\,  (\Sigma - {1\over2} \, {\rm I}) \, 
\, \big( \, \Delta t \,\, \Psi_1 (W)  + \Delta t^2 \,\, \Psi_2 (W) + \Delta t^3 \,\, \,\Psi_3 (W) \, \big) + {\rm O}(\Delta t^4) \, 
\end {array} \right. \monend
where $ \, {\rm I} \, $ is the identity matrix of dimension $ \, q - N $.
Observe that the two lines of (\ref{variables-microscopiques}) are equivalent thanks to (\ref{relaxation-Y})
and (\ref{Henon}). 
The vectors $\, \Phi_j(W) \ $ are analogous to $\, \Gamma_j (W)\, $ but not with the same dimension: $ \,  \Gamma_j(W) \in \RR^{N} \,$ 
and  $ \,  \Phi_j(W) \in \RR^{q-N} $. 

\newpage 
\smallskip  \monitem {\bf  Main result.   Recurrence formulas for fourth order expansion}

\noindent
We prove in this contribution that we have 
 the following recurrence formulas for the explicitation of the vectors $\, \Gamma_j (W) \, $ and $\, \Phi_j (W) \, $ 
up to fourth order accuracy:  
\moneq \label{formules} \left \{ \begin {array}{rl}
\Gamma_1 (W)  & \!\!\! =     A \, W + B \, \Phi(W) \\
\Psi_1 (W) & \!\!\! =   \dd \Phi(W) .  \Gamma_1 (W)  - \big( C \, W + D \, \Phi(W)  \big) \\ 
\Gamma_2 (W) & \!\!\! =  B \, \Sigma \, \Psi_1 (W) \\
\Psi_2 (W) & \!\!\! =   \Sigma \,  \dd \Psi_1 (W) .   \Gamma_1 (W) + \dd \Phi(W) .  \Gamma_2 (W) 
- D \, \Sigma \,  \Psi_1 (W)   \\ 
\Gamma_3(W)   & \!\!\! = B \, \Sigma  \, \Psi_2 (W) + {{1}\over{12}}  B_2 \, \Psi_1  (W)
-  {{1}\over{6}} \, B \,  \dd \Psi_1 (W) .  \Gamma_1 (W) \\
\Psi_3 (W) & \!\!\! = \Sigma \, \dd \Psi_1 (W) .  \Gamma_2 (W)  +  \dd \Phi(W) .  \Gamma_3(W) -  D \, \Sigma \, \Psi_2 (W) + \Sigma \, \dd \Psi_2 (W) .  \Gamma_1 (W) \\
&  +{1\over6} \, D \, \dd \Psi_1 (W) .  \Gamma_1 (W)  - {1\over12} \, D_2 \, \Psi_1 (W)  - {1\over12} \, \partial^2 \Psi_1 (W) .  \Gamma_1 (W)  \\ 
\Gamma_4(W)   & \!\!\! =      B \, \Sigma \, \Psi_3 (W) + {1\over4} \, B_2 \, \Psi_2 (W)
+  {1\over6} \, B \, D_2 \, \Sigma \, \Psi_1 (W)   -   {1\over6} \, A \, B \, \Psi_2 (W) \\
 &   -  {1\over6} \, B \, \big( \dd \gamma_1 (W) . \Gamma_2 (W) + \dd \gamma_2 (W) . \Gamma_1 (W) \big) 
 - {1\over6} \, B \, \Sigma \,  \partial^2 \Psi_1 (W) .  \Gamma_1 (W) \, . 
\end {array} \right. \monend
Observe that the order of the relations in (\ref{formules}) is mandatory. Each step has to be explicited
before the evaluation of the next expression. 
At second order accuracy, the three first lines of the  relations (\ref{formules}) 
summarizes  what is needed for a lot of applications
with only first and second order partial differential equations. 
Observe that the  compact form $ \,  \Gamma_2 = B \, \Sigma \, \Psi_1 \, $ for the second order term
put in evidence the important choice of a  H\'enon matrix  ``as  small as possible''. 
%
Observe that the relations  (\ref{formules}) have similarities with the ones obtained
in (\ref{formules-temps-continu}) for the  Chapman-Enskog expansion. But essentially 
the matrix $ \, \Sigma \, $ has replaced the matrix~$ \, Z^{-1} \, $ and new terms appear
due to the post-processing with Taylor expansions.

\bigskip \noindent {\bf Proposition 3. \quad  Equilibrium state and zero order expansion }  

\noindent
When $ \,  \Delta t \, $ tends to zero, the microscopic moments are close to their equilibrium value:
\moneq \label{y-ystar-ordre-0} 
Y = \Phi(W) + {\rm O}(\Delta t) \,,\quad Y^* = \Phi(W) + {\rm O}(\Delta t) \, .
\monend 

\monitem Proof of Proposition 3.

\noindent
We expand  one iteration of the scheme  (\ref{exponentielle-Lambda}) at  order zero:
\moneq  \label{m-m-star-ordre-0} 
m  + {\rm O}(\Delta t) = m^*  + {\rm O}(\Delta t^3) \, . 
\monend
Then we deduce from (\ref{m-m-star-ordre-0}): 
\moneqstar 
Y - Y^* =  {\rm O}(\Delta t) \, . 
\monendstar 
Then due to the iteration  (\ref{relaxation-Y}) and the fact that the matrix $ \, S \, $ is supposed fixed, 
we have  
\moneq \label{dvt-Y-ordre-0} 
Y =  \Phi(W)  + {\rm O}(\Delta t) \, . 
\monend
The first relation  of (\ref{y-ystar-ordre-0}) is satisfied. 
It is then immediate to deduce from the scheme iteration (\ref{relaxation-Y}) the second relation of (\ref{y-ystar-ordre-0}). 
The proposition is established. \hfill  $ \square $ 

\smallskip  \noindent 
In the following, we detail all the ingredients that conduct to the main result 
(\ref{edp-ordre-4}) to (\ref{formules}). 

\bigskip \bigskip    \noindent {\bf \large  5) \quad  Taylor expansion method at first and second order accuracy}    

\smallskip  \noindent 
In this section, we prove the two first orders of the expansion 
(\ref{edp-ordre-4})(\ref{variables-microscopiques})(\ref{formules}). 
We make the link with the Taylor expansion as presented in  \cite{Du07,Du08} and  
we apply the result for the isothermal D2Q9 scheme.

\newpage 
\bigskip \noindent {\bf Proposition 4. \quad  First order expansion }  
  
\noindent
When $ \,  \Delta t \, $ tends to zero, the macroscopic moments $ \, W \, $ satisfy asymptotically
the following system of $ \, N \, $ first order equations:
\moneq \label{edp-ordre-1} 
\partial_t W + \Gamma_1 (W)  = {\rm O}(\Delta t) \, . 
\monend
Moreover, the vector $ \,  \Gamma_1 (W) \, $ for the first order dynamics satisfy 
\moneq \label{Gamma-1} 
 \Gamma_1 (W)  = A \, W + B \, \Phi(W) \, . 
\monend
%

\monitem Proof of Proposition 4. 

\noindent
We expand  one iteration of the scheme  (\ref{exponentielle-Lambda}) at first order:
\moneq  \label{m-m-star-ordre-1}  
m + \Delta t \, \partial_t m + {\rm O}(\Delta t^2) = m^* - \Delta t \, \Lambda \, m^*  + {\rm O}(\Delta t^2) \, . 
\monend 
We can replace the vector $\, m \,$ by its two components $ \, W \, $ and $ \, Y $. We have for the first
line  
\moneqstar 
W +  \Delta t \, \partial_t W +  {\rm O}(\Delta t^2) =  \, W - \Delta t \, (A \, W + B \, Y^*)  + {\rm O}(\Delta t^2)  \, . 
\monendstar 
We simplify by the constant $ \, W \, $ term, 
divide by $ \, \Delta t \, $ and take into account the expansion~(\ref{y-ystar-ordre-0}): 
\moneqstar  
 \partial_t W + {\rm O}(\Delta t) =   - \big( A \, W + B \, \Phi(W) \big)  + {\rm O}(\Delta t) \, . 
\monendstar 
This relation is exactly the expansion (\ref{edp-ordre-1}) with the first order dynamics $ \, \Gamma_1 (W) \, $ 
evaluated according to the relation (\ref{Gamma-1}).

\smallskip  \monitem  {\bf Example of the D2Q9 fluid scheme}

The results presented here are essentially a reformulation 
of the classic work of Lallemand and Luo \cite{LL00}. 
With the moments introduced in (\ref{M-d2q9}), (\ref{d2q9-W}) and (\ref{d2q9-Y}), 
we denote by $ \,  \Phi_{\varepsilon} $, $ \,  \Phi_{xx}  \, $ and  $ \,  \Phi_{xy}  \, $
the equilibium values of the three first non conserved moments (\ref{d2q9-Y}). 
The dynamics at first order is given by (\ref{Gamma-1}). We have in this D2Q9 case, 
\moneqstar 
\Gamma_1 = \begin{pmatrix} \partial_x J_x + \partial_y J_y \\   \vspace{-.4 cm} \\ 
{2\over3} \, \lambda^2 \,  \partial_x \rho + {1\over6} \,  \partial_x  \Phi_{\varepsilon} +  {1\over2} \,  \partial_x  \Phi_{xx} 
+  \partial_y  \Phi_{xy} \\ \vspace{-.4 cm} \\ 
{2\over3} \, \lambda^2 \,  \partial_y \rho +{1\over6} \,  \partial_y  \Phi_{\varepsilon}  -  {1\over2} \,  \partial_y  \Phi_{xx}  
+  \partial_x  \Phi_{xy}  \end{pmatrix} \,. 
\monendstar 
These relations can be easily fitted with the Euler equations of gas dynamics in two space dimensions 
\moneqstar   \left \{ \begin {array}{l} 
 \partial_t  \rho +  \partial_x J_x + \partial_y J_y  \,=\, 0   \\ 
 \partial_t  J_x \,+\,  \partial_x \big(   {{1}\over{\rho}} \, J_x^2 + p \big)
 \,+ \,\partial_y \big(   {{1}\over{\rho}} \,  J_x \, J_y \big)  \,=\, 0  \\
\partial_t  J_y 
 \,+ \,\partial_x \big(  {{1}\over{\rho}} \,  J_x \, J_y \big)
 \,+\,  \partial_y \big(   {{1}\over{\rho}} \, J_y^2 + p \big)  \,=\, 0 \, . 
 \end {array} \right. \monendstar 

\smallskip \smallskip We identify the two expressions with $ \, J_x \equiv \rho \, u \, $ and $ \, J_y \equiv  \rho \, v$:  
\moneqstar   
{2\over3} \, \lambda^2 \, \rho \,+\,  {1\over6} \, \Phi_{\varepsilon} 
\,+\,  {1\over2} \, \Phi_{xx}  =  \rho \, u^2 + p  \,, \quad  
\Phi_{xy} =  \rho \, u \, v   \\ 
{2\over3} \, \lambda^2 \, \rho \,+\,  {1\over6} \, \Phi_{\varepsilon} 
\,-\,  {1\over2} \, \Phi_{xx}  =  \rho \, v^2 + p  \,, \quad 
\monendstar 
and we deduce the expression for the three first noncenserved moments: 
\moneqstar 
 \Phi_{\varepsilon}  =  6 \, p -  4 \, \lambda^2 \, \rho \, + \, 3 \, \rho \, ( u^2 + v^2 )     \,, \quad  
\Phi_{xx}  =   \rho \, ( u^2 - v^2 )     \,, \quad  
 \Phi_{xy} = \rho \, u \, v    \, . 
\monendstar 
The other nonequilibrium moments $ \,  q_x $, $ \,  q_y  \, $ and $ \, h \, $ 
does not play any role in the first order partial equivalent equations. 

\newpage 
\bigskip \noindent {\bf Proposition 5. \quad  Taylor expansion method at second order accuracy } 

\noindent
An essential result for the applications is the second order 
asymptotic analysis for a time step $ \,  \Delta t \, $ (or a space step $ \, \Delta x$) tending to zero.
The expansion of the microscopic variables takes the form 
\moneq \label{variables-microscopiques-ordre-1} 
 Y = \Phi(W) +  \Delta t \,\,  S^{-1} \, \Psi_1 (W)  + {\rm O}(\Delta t^2) 
\monend
and the vector $ \,  \Psi_1 (W) \, $ satisfies the relation
\moneq \label{Phi-1} 
\Psi_1 (W)  = \dd \Phi(W) \, . \,  \Gamma_1 (W)  - \big( C \, W + D \, \Phi(W)  \big) \, . 
\monend
A set of second order partial differential equations emerges 
\moneq \label{edp-ordre-2} 
\partial_t W + \Gamma_1 (W) + \Delta t \,\, \Gamma_2 (W)  = {\rm O}(\Delta t^2) \, . 
\monend
The vector  $ \, \Gamma_1 (W) \, $ has been precised in (\ref{Gamma-1}) 
the vector  $ \, \Gamma_2 (W) \, $ is  obtained after
two  derivations  of the conserved moments $ \, W \, $ and 
the equilibrium vector $ \, \Phi(W) $:
\moneq \label{Gamma-2} 
\Gamma_2 (W)  = B \,\, \Sigma \,\, \Psi_1  (W) \, . 
\monend
The  expression of the second order term puts in evidence the importance of taking over-relaxation in the applications
where $ \, \Sigma \, $  should also be as small as possible.  

\monitem Proof of Proposition 5.

\noindent 
We first consider the first order expansion of microscopic moments.  
From the expansion (\ref{m-m-star-ordre-1}) of one iteration of the scheme  at first order,
we can extract relations for the second component.
We have for the microscopic moments 
\moneqstar 
 Y +  \Delta t \,\, \partial_t Y  + {\rm O}(\Delta t^2) =  Y^*  - \Delta t \,\, (C \, W + D \, Y^*)  + {\rm O}(\Delta t^2) \, . 
\monendstar 
Then 
\moneqstar  
Y - Y^*  = -  \Delta t \, \Big( \partial_t Y +  \big(C \, W + D \, \Phi(W) \big) \Big) + {\rm O}(\Delta t^2) \, . 
\monendstar 
The explicitation of $ \,  \partial_t Y \, $ is given by the chain rule from the expansion (\ref{dvt-Y-ordre-0}): 

\smallskip \noindent $ \displaystyle 
 \partial_t Y =  \partial_t \big( \Phi(W) +  {\rm O}(\Delta t) \big) $

\smallskip \noindent $ \qquad \displaystyle =  
\dd \Phi(W)  \, . \,  \partial_t W +  {\rm O}(\Delta t)   $

\smallskip \noindent $ \qquad \displaystyle =  
 \dd \Phi(W) .  (-\Gamma_1 )   +  {\rm O}(\Delta t) \, . $

\smallskip \noindent 
Then 
\moneqstar 
 S \,\,  \big( Y - \Phi(W) \big) = \Delta t \,\, \big( \dd \Phi(W) \, . \, \Gamma_1 - (C \, W + D \, \Phi) \big)  + {\rm O}(\Delta t^2)  \, . 
\monendstar 
The relation   (\ref{variables-microscopiques-ordre-1}) is established, and the expression of $ \, \Psi_1 (W) \, $ is 
obtained by the relation~(\ref{Phi-1}).
Before using the expansion (\ref{variables-microscopiques-ordre-1}), we must expand the microscopic moments $ \, Y^* \, $
after relaxation up to first order accuracy. We have the following calculus 

 \smallskip  \noindent $ \displaystyle 
Y^* =  Y + S \, ( \Phi(W)  - Y)  $ 

 \smallskip  \noindent \quad $ \,\,\, \displaystyle = 
Y - \big( \Delta t \, \Psi_1 +  {\rm O}(\Delta t^2)  \big) $ \hfill due to (\ref{variables-microscopiques-ordre-1}) 

 \smallskip  \noindent \quad $ \,\,\, \displaystyle = 
\Phi(W) + \Big( \Sigma + {{1}\over{2}} \, {\rm I} \Big)  \, \big( \Delta t \, \Psi_1 (W) +  {\rm O}(\Delta t^2)  \big) 
 -  \big( \Delta t \, \Psi_1 +  {\rm O}(\Delta t^2)  \big)   $ 
 \hfill due to (\ref{Henon}) and (\ref{variables-microscopiques-ordre-1}) 

 \smallskip  \noindent \quad $ \,\,\, \displaystyle = 
\Phi(W) + \Big( \Sigma - {{1}\over{2}} \, {\rm I} \Big)   \,  \Delta t \, \Psi_1(W)   + {\rm O}(\Delta t^2)  $

\smallskip  \noindent
and we have the relation 
\moneq \label{Y-star-ordre-1}  
Y^* =  \Phi(W) + \Big( \Sigma - {{1}\over{2}} \, {\rm I} \Big)  \,   \Delta t \,\,    \Psi_1(W)  + {\rm O}(\Delta t^2) \, . 
\monend
The two relations  (\ref{variables-microscopiques-ordre-1}) and   (\ref{Y-star-ordre-1})
are nothing else that the first order terms
of the general expansion (\ref{variables-microscopiques}).

\smallskip \monitem  Second order partial differential equations 

\noindent 
We expand now one iteration of the scheme  (\ref{exponentielle-Lambda}) at second order accuracy:
\moneq     \label{m-m-star-ordre-2}   
m + \Delta t \, \partial_t m + {1\over2} \, \Delta t^2 \, \partial_t^2 m + {\rm O}(\Delta t^3) = 
m^* - \Delta t \, \Lambda \, m^* + {1\over2} \, \Delta t^2 \, \Lambda^2 \,  m^* + {\rm O}(\Delta t^3) \, . 
\monend 
We replace the vector $\, m \,$ by its two components $ \, W \, $ and $ \, Y $.
We use the decomposition (\ref{bloc-Lambda})  of the momentum-velocity operator matrix $ \, \Lambda $. 
We have for the first component 
\moneq \label{dvt-W-brut}   \left \{ \begin {array}{l}  \displaystyle 
W +  \Delta t \, \partial_t W + {1\over2} \, \Delta t^2 \, \partial_t^2 W + {\rm O}(\Delta t^3) = \\    \displaystyle 
\qquad \qquad  \, W - \Delta t \, (A \, W + B \, Y^*) + {1\over2} \, \Delta t^2 \, (A_2 \, W + B_2 \, Y^*)   + {\rm O}(\Delta t^3) \, . 
\end {array}  \right. \monend
We  explicit some terms present in the relation (\ref{dvt-W-brut}). 
We have 
\moneqstar 
\partial_t W  =  -\Gamma_1 (W) +  {\rm O}(\Delta t)  = -(A \, W + B \,  \Phi(W) ) +  {\rm O}(\Delta t) 
\monendstar 
Then, using the formal rule stated at the beginning of Section~4,  
we keep the order of accuracy after derivation: 

\smallskip \noindent $ \displaystyle 
\partial^2_t W = -\partial_t \big( \Gamma_1 +  {\rm O}(\Delta t) \big) $ 

\smallskip \noindent $ \qquad \,\, \displaystyle = 
 -\dd \Gamma_1 \, . \,  \partial_t W +   {\rm O}(\Delta t)   $ 

\smallskip \noindent $ \qquad \,\, \displaystyle = 
\dd \Gamma_1 \, . \, \Gamma_1  +   {\rm O}(\Delta t)  $ 

\smallskip \noindent $ \qquad \,\, \displaystyle = 
A \,  \Gamma_1  + B \, \dd \Phi \,.\,  \Gamma_1    +   {\rm O}(\Delta t)  \, . $

\smallskip  \noindent 
We deduce from (\ref{dvt-W-brut}): 
\moneqstar   
\partial_t W = - {1\over2} \, \Delta t \, \partial_t^2 W  - ( A \, W + B \, Y^* ) 
+ {1\over2} \, \Delta t \, (A_2 \, W + B_2 \, Y^*)   + {\rm O}(\Delta t^2) \, . 
\monendstar 
We order the various terms by powers of $ \, \Delta t $:

\smallskip \noindent $ \displaystyle 
\partial_t W = -A \, W - B \, \Big( \Phi + (\Sigma - {1\over2} {\rm I} )\, \Delta t \, \Psi_1 \Big) 
- {1\over2} \, \Delta t \, (  A \,  \Gamma_1  + B \, \dd \Phi .  \Gamma_1 ) $

\smallskip \noindent $ \qquad \qquad  \qquad  \displaystyle 
  + {1\over2} \, \Delta t \, \big( (A^2 + B \, C) \, W + (A \, B + B \, D) \, \Phi  \big)  + {\rm O}(\Delta t^2) $

\smallskip \noindent $ \qquad \, \displaystyle =   
-A \, W - B \, \Phi  +  \Delta t \, \Big[  -B \, \Sigma \, \Psi_1 
+ {1\over2} B \, \big(  \dd \Phi \,.\,  \Gamma_1 - C \, W - D \, \Phi \big) 
- {1\over2}  A \,\, (  A \, W +  B \,  \Phi)  $

\smallskip \noindent $ \qquad \qquad \qquad \displaystyle 
 -   {1\over2} \, B \,  \dd \Phi .  \Gamma_1 
+  {1\over2} \,  (A^2 + B \, C) \, W +  {1\over2} \,  (A \, B + B \, D) \, \Phi   \Big]  + {\rm O}(\Delta t^2) $ 

\smallskip \noindent $ \qquad \, \displaystyle =   
 - \Gamma_1 (W)  -  \Delta t \, B \,\, \Sigma \,\, \Psi_1 (W)  + {\rm O}(\Delta t^2) $ 

\smallskip \noindent 
and the relation (\ref{edp-ordre-2}) is proven with   $ \, \Gamma_2 (W)  =  B \,\, \Sigma \,\, \Psi_1 \, $ as 
proposed in  (\ref{Gamma-2}). \hfill $\square$

\smallskip \monitem  {\bf Link with the original Taylor expansion method }

\noindent 
In \cite{Du07,Du08}, the moments at equilibrium are denoted  $ \, m_\ell^{\rm eq}  $:  
\moneqstar
m_\ell^{\rm eq} (W)  = \left \{ \begin {array}{l}
\displaystyle W_\ell \qquad \qquad  \quad    {\rm if \, } \, 0 \leq \ell < N \\ 
\displaystyle \big(\Phi (W) \big)_{\ell-N}   \quad \,\, {\rm if \, } \, \ell \geq N \,, 
\end{array}  \right. \monendstar 
the  coefficients $ \,  \widetilde{\Lambda}_{i \ell}^\beta \, $ are given according to 
\moneq  \label{Lambda-2007} 
\widetilde{\Lambda}_{i \ell}^\beta = \sum_j  M_{ij} \,  v_j^\beta \, (M^{-1})_{j \ell} \,, \quad 1 \leq \beta \leq d \,, 
0 \leq i,\, \ell < q \,, 
\monend 
are closely related to the operator matrix $ \, \Lambda \, $ introduced in (\ref{Lambda}). 
The coefficients of the H\'enon's matrix (\ref{Henon}) have a shifted numbering 
\moneq   \label{sigma-henon} 
\widetilde{\sigma}_\ell = \sigma_{\ell - N} \equiv {{1}\over{s_{\ell-N}}} - {1\over2} \,,\quad \ell \geq N 
\monend
and the defect of conservation $ \,\widetilde{\theta}_k \, $  is defined according to 
\moneq \label{defaut-conservation} 
\widetilde{\theta}_k (W) \equiv \partial_t m_k^{\rm eq} + \sum_{0 \leq \ell < q} \,
\sum_{1 \leq \beta \leq d}  \widetilde{\Lambda}_{k \ell}^\beta \,\, \, \partial_\beta  m_\ell^{\rm eq}(W)  \,,\quad k \geq N \, . 
\monend
%

\bigskip \noindent {\bf Proposition 6. \quad  Equivalence at second order accuracy }

\noindent 
The second order expansion (\ref{edp-ordre-2}), 
(\ref{variables-microscopiques-ordre-1}),  (\ref{Gamma-1}), (\ref{Phi-1}) and (\ref{Gamma-2})
is equivalent to  the expansion proposed in \cite{Du07,Du08}:
\moneq \label{edp-ordre-2-2007} 
  \partial_t W_i + \widetilde{\Lambda}_{i \ell}^\beta \,\, \partial_\beta m_\ell^{\rm eq} (W) 
=  \Delta t \, \sum_{k \geq N}   \widetilde{\Lambda}_{i k}^\beta \, \widetilde{\sigma}_k \,\,\partial_\beta \widetilde{\theta}_k (W) 
+  {\rm O}(\Delta t^2)  \,. 
\monend

\monitem Proof of Proposition 6.

\noindent 
Remark first that
\moneqstar 
\Lambda_{i k} =    \sum_{j=0}^{q-1}  M_{ij} \,  \Big( \sum_{1 \leq \beta \leq d} v_j^\beta \, \partial_\beta \Big) \, (M^{-1})_{j k}
\,=\,\sum_{1 \leq \beta \leq d} \, \Big(   \sum_{j=0}^{q-1}  M_{ij} \,  v_j^\beta \,  (M^{-1})_{j k} \Big) \,  \partial_\beta 
\,=\,\sum_{1 \leq \beta \leq d} \widetilde{\Lambda}_{i k}^\beta \,  \partial_\beta \, . 
\monendstar 
Then on one hand, the component number $ \, i \, $ of the first order term
$\,\big(  A \, W + B \, \Phi(W) \big) \, $ can be expanded in the following way:

\smallskip \noindent $ \displaystyle
\big(  A \, W + B \, \Phi(W) \big)_i = \sum_{k < N}\widetilde{\Lambda}_{i k}^\beta \,\,  \partial_\beta W_k +
\sum_{k \geq N}\widetilde{\Lambda}_{i k}^\beta \,\,  \partial_\beta \Phi_k
$ 

\smallskip \noindent $ \displaystyle \qquad \qquad  \qquad \qquad \, = 
\sum_{k < N}\widetilde{\Lambda}_{i k}^\beta \,\,  \partial_\beta  m_k^{\rm eq}  +
\sum_{k \geq N}\widetilde{\Lambda}_{i k}^\beta \,\,  \partial_\beta m_k^{\rm eq}  $

\smallskip \noindent $ \displaystyle \qquad \qquad \qquad  \qquad \, = 
\sum_{0 \leq k < q}\widetilde{\Lambda}_{i k}^\beta \,\,  \partial_\beta  m_k^{\rm eq}  \,  $

\smallskip \noindent
and the first order terms of  (\ref{edp-ordre-2}) and (\ref{edp-ordre-2-2007}) are identical.
Secondly, due to (\ref{defaut-conservation}), we have for $ \, k \geq N \, $ 

\smallskip \noindent $ \displaystyle
\widetilde{\theta}_k =  \partial_t m_k^{\rm eq} + \sum_{\ell \beta} \widetilde{\Lambda}_{k \ell}^\beta \,\,  \partial_\beta  m_\ell^{\rm eq}  $

\smallskip \noindent $ \displaystyle \quad \, = \partial_t \Phi_{k-N}  + \big( \Lambda . m^{\rm eq} \big)_k  $ 

\smallskip \noindent $ \displaystyle \quad \, =
\big( \dd \Phi \, . \, \partial_t W \big)_{k-N}  + \big( \Lambda . m^{\rm eq} \big)_k  $ 

\smallskip \noindent $ \displaystyle \quad \, =
\big( \dd \Phi \, . \, ( \Gamma_1 + {\rm O} (\Delta t))  \big)_{k-N}   - \big( C \, W + D \, \Phi \big)_{k-N}  $

\smallskip \noindent
and due to (\ref{Phi-1}),  
\moneqstar
\widetilde{\theta}_k = - \big( \Psi_1  \big)_{k-N}  + {\rm O}(\Delta t) \,, \quad k \geq N \, . 
\monendstar 
We deduce that for second order terms with $ \, 0 \leq i < N \, $: 

\smallskip \noindent $ \displaystyle
\big(  B \,\, \widetilde{\sigma} \,\, \Psi_1 \big)_i = \sum_{0 \leq \ell < q-N} B_{i \ell} \, \widetilde{\sigma}_{\ell+N} \, \, \big( \Psi_1(W) \big)_\ell $ 

\smallskip \noindent $ \displaystyle \qquad  \qquad \,\,\,\,   =
\sum_{k \, \geq N} \, \sum_{0 \leq \beta \leq d} \widetilde{\Lambda}_{i k}^\beta \,\, \partial_\beta
\Big( \sigma_k \,  \big(-\widetilde{\theta}_k +  {\rm O}(\Delta t) \big) \Big) $ \hfill with $ \, k = \ell + N $ 

\smallskip \noindent $ \displaystyle \qquad  \qquad \,\,\,\,   =
- \sum_{k \, \geq N} \, \sum_{0 \leq \beta \leq d} \widetilde{\Lambda}_{i k}^\beta \, 
\widetilde{\sigma}_k  \, \big( \partial_\beta \widetilde{\theta}_k  \big)  +  {\rm O}(\Delta t) $ 

\smallskip \noindent
because the coefficients $ \, \widetilde{\sigma}_k \, $ are supposed constant.
Then the relation (\ref{edp-ordre-2-2007}) is established.  \hfill $\square$

\smallskip  \monitem  {\bf Fluid D2Q9  diffusive tensor}

\noindent 
If we compute the second order term $ \, \Gamma_2 \, $ for the D2Q9 scheme introduced previously, 
the holy grail  would be to recover the viscous terms of the compressible Navier Stokes equations in two space
dimensions.
More precisely,  
with a given shear viscosity $ \, \mu \, $ and a given bulk viscosity $\, \zeta $, 
we write the viscous terms of the lattice Boltzmann expansion as 
\moneq   \label{erreur-NS-d2q9} 
- \Delta t \, \Gamma_2 = \begin {pmatrix} 0  \\ \partial_j \tau_{xj} 
 \equiv \partial_x ( 2 \, \mu \, \partial_x u 
+ (\zeta-\mu)  ( \partial_x u   +  \partial_y v  ) ) + \partial_y (\mu ( \partial_x v   +  \partial_y u ) ) 
\\ \partial_j \tau_{yj}  \equiv \partial_x (\mu ( \partial_x v  +  \partial_y u  ) ) 
+ \partial_y (  (\zeta-\mu) ( \partial_x u   +  \partial_y v  ) + 2 \, \mu \,  \partial_y v  ) 
\end{pmatrix} +  \begin {pmatrix} 0  \\ \Psi_x \\ \Psi_y \end{pmatrix} \, ,  
\monend
with the notation $ \, \Psi_x $,  $ \, \Psi_y \, $ for the 
error for the dissipation  of momentum.  

\smallskip \noindent 
The problem to enforce $ \, \Psi_x = \Psi_y = 0 \, $ is ill-posed, as recognized by 
Dellar \cite{De02, De14}, H\'azi and K\'avr\'an \cite {HK06} and 
Prasianakis {\it at al.} \cite{PKMB09} among others. 
Nevertheless, it can be proven that the unphysical terms 
involving second order space derivatives involving  the density 
can be removed with a pressure law of the type 
\moneqstar 
p(\rho) = {{\lambda^2}\over{3}} \, \rho \,, 
\monendstar  
{\it id est} a sound velocity equal to $ \, {{1}\over{\sqrt{3}}} \, $ in the velocity unit of the lattice Boltzmann scheme. 
Nevertheless, with the notation  $ \,  \Phi_{qx} \, $ and   $ \,  \Phi_{qy} \, $ for the equilibrium fonctions
of the moments $ \, q_x \, $ and $ \, q_y \, $ and 
the coefficients of the H\'enon matrix $\, \Sigma \, $ fixed as
\moneq \label{Sigma-d2q9} 
\Sigma = {\rm diag} \big( \sigma_e ,\, \sigma_x ,\, \sigma_x ,\, \sigma_q ,\, \sigma_q  ,\, \sigma_h  \big)  \, , 
\monend
the choice 
\moneqstar  \left \{ \begin {array}{l}
 \Phi_{qx} =  - \rho \, \lambda^2\, u  + 3 \, \rho \, (u^2 + v^2) \, u  \\ 
 \Phi_{qy} =  - \rho \, \lambda^2\, v  + 3 \, \rho \, (u^2 + v^2) \, v \,, 
\end{array} \right.  \monendstar   
conducts to the relation (\ref{erreur-NS-d2q9}) with associated viscosities fixed in a very classical way: 
\moneq  \label{viscosites-d2q9}  
\mu    =  {{\lambda}\over{3}} \,  \rho \, \sigma_x \, \Delta x  \,, \quad
\zeta    = {{\lambda}\over{3}} \,  \rho  \, \sigma_e \, \Delta x  \, .  
\monend
The equilibrium value $ \, \Phi_h \, $ of the last fourth order  moment $ \, h \, $ (see the relation (\ref{d2q9-Y}))
 has no incidence on the value of the diffusive term $ \, \Gamma_2 $. 
The error is then reduced to third order terms relative to the velocity and we have precisely 
\moneqstar  
\left \{ \begin {array}{l}
\Psi_x =  \sigma_x \, \Delta t \, \Big[ \partial_x \big( u^3 \, \partial_x \rho - v^3 \, \partial_y \rho 
+ 3 \, \rho \, (u^2 \, \partial_x u  - v^2 \, \partial_y v) \big) \\ \qquad  \qquad  \qquad  \qquad 
+ \partial_y \big( -v^3 \, \partial_x \rho - u^3 \,  \partial_y \rho 
-3 \, \rho \, ( u^2 \, \partial_y u + v^2 \, \partial_x v ) \big) \Big] \\
\Psi_y =  \sigma_x \, \Delta t \, \Big[ \partial_x \big( -v^3 \, \partial_x \rho 
- u^3 \, \partial_y \rho -3 \, \rho \, ( u^2 \, \partial_y u + v^2 \, \partial_x v \big) \\ \qquad  \qquad  \qquad  \qquad  
+ \partial_y \big(-u^3 \,  \partial_x \rho + v^3 \,  \partial_y \rho + 3 \, \rho \, (-u^2 \, \partial_x u + v^2 \,  \partial_y v )  \big) \Big] \,.
\end{array} \right.   \monendstar   

\newpage 
\bigskip \bigskip    \noindent {\bf \large  6) \quad  Taylor expansion method at third and fourth order accuracy}    

\noindent
The two  following orders of the  Taylor expansion 
does not set any theoretical difficulty, except the care of algebra with not so short formal expressions...

\bigskip \noindent {\bf Proposition 7. \quad   Nonlinear third-order expansion  } 

\noindent
When the time step $ \,  \Delta t \, $ has an infinitesimal value, 
the expansion of the microscopic variables takes the form 
\moneq \label{Y-ordre-2} 
 Y = \Phi(W) +  \Delta t \,\,  S^{-1} \, \Psi_1 (W) +  \Delta t^2 \,\,  S^{-1} \, \Psi_2 (W)    + {\rm O}(\Delta t^3) \, . 
\monend
The vector $ \,  \Psi_1 (W) \, $ is evaluated at the relation (\ref{Phi-1})
and $ \,  \Psi_2 (W) \, $   satisfies the relation
\moneq \label{Phi-2} 
\Psi_2  (W)  = \Sigma \,  \dd \Psi_1 (W) . \Gamma_1 (W) + \dd \Phi(W) .  \Gamma_2 (W) - D \, \Sigma \,  \Psi_1 (W)  \, . 
\monend
It is possible to precise a system of third order partial differential equations  
\moneq \label{edp-ordre-3} 
\partial_t W + \Gamma_1 (W) + \Delta t \, \Gamma_2 (W)  + \Delta t^2 \, \Gamma_3 (W)  = {\rm O}(\Delta t^3) \, . 
\monend
The vectors  $ \, \Gamma_1 (W) \, $ and  $ \, \Gamma_2 (W) \, $
have been precised at the relations  (\ref{Gamma-1}) and (\ref{Gamma-2}) respectively.
We have moreover 
\moneq \label{Gamma-3} 
\Gamma_3 (W)  = B \, \Sigma \, \Psi_2 (W) + {1\over12} \, B_2 \, \Psi_1 (W)
- {1\over6} \,B \, \dd \Psi_1 (W) .  \Gamma_1(W)  \, . 
\monend

\monitem Proof of Proposition 7.  

\noindent
We first establish the relations (\ref{Y-ordre-2}) and (\ref{Phi-2}). 
We start by the expansion of the scheme (\ref{m-m-star-ordre-2}) at order 2.
For the second component of nonconserved moments, we have

\smallskip \noindent  $ 
Y +  \Delta t \, \partial_t Y + {1\over2} \, \Delta t^2 \, \partial_t^2 Y + {\rm O}(\Delta t^3) =
Y^*  - \Delta t \, (C \, W + D \, Y^*) + {1\over2} \, \Delta t^2 \, (C_2 \, W + D_2 \, Y^*)   + {\rm O}(\Delta t^3) $

\smallskip \noindent 
and, due to the relaxation process $ \, Y - Y^*  \equiv  S \, ( Y - \Phi (W) ) $, we have 
\moneq    \label{dvt-Y-brut-2}
S \, ( Y - \Phi (W) )  =
- \Delta t \, ( C \, W + D \, Y^* + \partial_t Y) 
+  {1\over2} \, \Delta t^2 \, (C_2 \, W + D_2 \, Y^* - \partial_t^2 Y ) + {\rm O}(\Delta t^3) \, . \monend 
We can precise the two first time derivatives of the non-conserved moments. 
The relation (\ref{variables-microscopiques-ordre-1}) can also be written as
\moneqstar 
Y = \Phi(W) + \Big( \Sigma + {{1}\over{2}} \, {\rm I} \Big)   \,  \Delta t \, \Psi_1 (W)  + {\rm O}(\Delta t^2) \, . 
\monendstar  
We differentiate this previous relation relative to time:   

\smallskip \noindent  
 $ \, \partial_t Y = \dd \Phi . \, \partial_t W + 
\big( \Sigma + {{1}\over{2}} \, {\rm I} \big)   \,\,   \Delta t \,\,   \dd \Psi_1 . \, \partial_t W + {\rm O}(\Delta t^2)  $

\smallskip \noindent  \qquad $ = 
- \dd \Phi . \,  \big(  \Gamma_1 (W) + \Delta t \, \Gamma_2 (W)  + {\rm O}(\Delta t^2)  \big) 
+ \Delta t \,  \big( \Sigma + {{1}\over{2}} \, {\rm I} \big)   \,  \dd \Psi_1 . \,  \big(  -\Gamma_1 +  {\rm O}(\Delta t)  \big)
+ {\rm O}(\Delta t^2)  $


and 
\moneq    \label{dtY-ordre-1}
\partial_t Y = - \dd \Phi(W)  . \Gamma_1 (W) - \Delta t \, \Big[
\dd \Phi(W)  . \Gamma_2 (W) + \Big( \Sigma + {{1}\over{2}} \, {\rm I} \Big) \,\,   \dd \Psi_1  (W) . \Gamma_1 (W) \Big] 
+ {\rm O}(\Delta t^2) \,. 
\monend
Because we need the second time derivative $\, \partial_t^2 Y \, $ in the right hand side of (\ref{dvt-Y-brut-2}), 
we differentiate the relation (\ref{dtY-ordre-1}) relative to time.
Using the notation (\ref{notation-gamma-j}), we deduce 
\moneq    \label{dt2Y-ordre-0}
\partial_t^2 Y =  \dd \gamma_1 . \, \Gamma_1   + {\rm O}(\Delta t) \, . 
\monend
We precise now the expansions of $ \,  C \, W + D \, Y^* \, $ and $ \, C_2 \, W + D_2 \, Y^* $. We have 

\smallskip \noindent  
$ C \, W + D \, Y^* =  C \, W + D  \, \big[ \Phi(W) + \big( \Sigma - {{1}\over{2}} \, {\rm I} \big)   \, 
 \Delta t \, \Psi_1  \big] + {\rm O}(\Delta t^2) $

\noindent  and  
\moneq    \label{CWDYst-ordre-1}
C \, W + D \, Y^* =  \gamma_1 (W) -  \Psi_1 (W) 
+ \Delta t \,\, D \, \Big( \Sigma - {{1}\over{2}} \, {\rm I} \Big) \, \Psi_1 (W)   +  {\rm O}(\Delta t^2)  \, . 
\monend
With a new set of operator matrices, 

\newpage 
\smallskip \noindent  
$ C_2 \, W + D_2 \, Y^* =  C_2 \, W + D_2  \big( \Phi(W) +  {\rm O}(\Delta t) \big)  $

\smallskip \noindent \qquad \qquad \qquad  $ \,\,\, = 
(C \, A + D \, C) \, W +  (C \, B + D^2) \,\Phi(W) +   {\rm O}(\Delta t)  $  

\smallskip \noindent \qquad \qquad \qquad  $ \,\,\, = 
C \, \Gamma_1 + D \, ( \dd \Phi . \Gamma_1 - \Psi_1 ) +   {\rm O}(\Delta t)  $  

\noindent  and  because $ \,\, \,\dd \Phi . \Gamma_1 = \dd \gamma_1 . \Gamma_1 - ( C \, \Gamma_1 + D \, \dd \Phi . \Gamma_1 ) $, 
\moneq    \label{C2WD2Yst-ordre-0}
C_2 \, W + D_2 \, Y^* =  \dd ( \gamma_1 - \Psi_1 ) . \Gamma_1 (W) -  D \, \Psi_1  (W)  +   {\rm O}(\Delta t) \, . 
\monend
Then due to
(\ref{dvt-Y-brut-2}),  (\ref{dtY-ordre-1}), (\ref{dt2Y-ordre-0}), (\ref{CWDYst-ordre-1})
and (\ref{C2WD2Yst-ordre-0}), we have

\smallskip \noindent  $
S \, ( Y - \Phi (W) )  =
- \Delta t \, ( C \, W + D \, Y^* + \partial_t Y) 
+  {1\over2} \, \Delta t^2 \, (C_2 \, W + D_2 \, Y^* - \partial_t^2 Y ) + {\rm O}(\Delta t^3) $ 

\smallskip \noindent  \quad $ =
- \Delta t \, \big( \gamma_1   -  \Psi_1   
+ \Delta t \,\, D \, \big( \Sigma - {{1}\over{2}} \, {\rm I} \big) \, \Psi_1   \big) 
- \, \Delta t \, \big( - \dd \Phi   . \Gamma_1   - \Delta t \,
\big[ \dd \Phi   . \Gamma_2   + \big( \Sigma + {{1}\over{2}} \, {\rm I} \big) \,\,   \dd \Psi_1 . \Gamma_1   \big] \big) $

\smallskip \noindent  \qquad $ 
+  {1\over2} \, \Delta t^2 \, \big(  \dd ( \gamma_1 - \Psi_1 ) . \Gamma_1   -  D \, \Psi_1    \big)
- {1\over2} \, \Delta t^2 \,\, \dd \gamma_1 . \, \Gamma_1  +   {\rm O}(\Delta t^3) $

\smallskip \noindent  \quad $ = 
\Delta t \,\Psi_1  + \Delta t^2 \, \big( - D \, \Sigma \, \Psi_1  +  {1\over2} \, D \, \Psi_1  +  \dd \Phi   . \Gamma_2 
+ \Sigma \,  \dd \Psi_1 . \Gamma_1  -  {1\over2} \, D \, \Psi_1  \big) +  {\rm O}(\Delta t^3) $

\smallskip \noindent  \quad $ = 
\Delta t \,\Psi_1  + \Delta t^2 \, \big(  - D \, \Sigma \, \Psi_1  +  \dd \Phi   . \Gamma_2  + \Sigma \,  \dd \Psi_1 . \Gamma_1 \big)
+  {\rm O}(\Delta t^3) $  

and the relations (\ref{Y-ordre-2}) and (\ref{Phi-2}) are established.

\smallskip  \monitem
Third order partial differential equations  

\noindent
We consider now the Taylor expansion of the relation (\ref{exponentielle-Lambda}) at third order: 
\moneq     \label{m-m-star-ordre-3}   \left \{ \begin {array}{l} 
m + \Delta t \, \partial_t m + {1\over2} \, \Delta t^2 \, \partial_t^2 m + {1\over6} \, \Delta t^3 \, \partial_t^3 m  + {\rm O}(\Delta t^4) 
= \\  \qquad   \qquad    m^* - \Delta t \, \Lambda \, m^* + {1\over2} \, \Delta t^2 \, \Lambda^2 \,  m^*   
 -{1\over6} \, \Delta t^3 \, \Lambda^3 \,  m^* +  {\rm O}(\Delta t^4)  \, . 
\end{array} \right. \monend 
We consider the first conserved components of the moments:
\moneqstar      \left \{ \begin {array}{l} 
W +  \Delta t \, \partial_t W + {1\over2} \, \Delta t^2 \, \partial_t^2 W 
+ {1\over6} \, \Delta t^3 \, \partial_t^3 W + {\rm O}(\Delta t^4) =  W - \Delta t \, (A \, W + B \, Y^*) 
 \\  \qquad   \qquad    + {1\over2} \, \Delta t^2 \, (A_2 \, W + B_2 \, Y^*)
-  {1\over6} \, \Delta t^3 \, (A_3 \, W + B_3 \, Y^*)   + {\rm O}(\Delta t^4) \,. 
\end{array} \right. \monendstar 
We simplify the constant term $ \, W \, $ and divide by $ \, \Delta t $. We deduce
\moneq \label{dvt-W-brut-ordre-3}   \left \{ \begin {array}{l}  
\partial_t W =  -A \, W - B \, Y^* +  {1\over2} \, \Delta t \, (A_2 \, W + B_2 \, Y^* - \partial_t^2 W )
\\  \qquad   \qquad - {1\over6} \, \Delta t^2 \, (A_3 \, W + B_3 \, Y^* + \partial_t^3 W ) + {\rm O}(\Delta t^3)  \, . 
\end{array} \right. \monend
We explicit the partial derivatives $ \, \partial_t^2 W  \, $ and $ \, \partial_t^3 W \, $ at orders one
and zero respectively. From the relation 
\moneqstar
\partial_t W + \Gamma_1 + \Delta t \, \Gamma_2  = {\rm O}(\Delta t^2)  
\monendstar 
we have

\smallskip \noindent  $
\partial_t^2 W = \partial_t \, \big( -\Gamma_1 - \Delta t \, \Gamma_2  + {\rm O}(\Delta t^2)  \big)   $ 

\smallskip \noindent  \qquad  $  \,\, = 
\dd \big( -\Gamma_1 - \Delta t \, \Gamma_2 + {\rm O}(\Delta t^2) \big) . (\partial_t W) $ 

\smallskip \noindent  \qquad  $  \,\, =
\dd \big( -\Gamma_1 - \Delta t \, \Gamma_2 + {\rm O}(\Delta t^2) \big) \, . \, \big( -\Gamma_1 - \Delta t \, \Gamma_2
+ {\rm O}(\Delta t^2)  $

\smallskip \noindent
and 
\moneq    \label{dt2W-ordre-1}
\partial_t^2 W = 
\dd  \Gamma_1 .\, \Gamma_1 + \Delta t \, ( \dd  \Gamma_1 .\, \Gamma_2 + \dd  \Gamma_2 .\, \Gamma_1 ) + {\rm O}(\Delta t^2) \, . 
\monend
Then 

\smallskip 
$\, \partial_t^3 W = \partial_t (  \partial_t^2 W )  $

\smallskip  \qquad 
$  \,\,\, =  \partial_t \big( \dd  \Gamma_1 .\, \Gamma_1  \big) + {\rm O}(\Delta t) $  

\smallskip  \qquad 
$ \,\,\, = \dd \big( \dd  \Gamma_1 .\, \Gamma_1 \big) .  \partial_t W  + {\rm O}(\Delta t) $  

\smallskip  \qquad 
$ \,\,\, = -\dd \big( \dd  \Gamma_1 .\, \Gamma_1 \big) . \Gamma_1 + {\rm O}(\Delta t) $  
\moneq    \label{dt3W-ordre-0}
\partial_t^3 W = -\partial^2 \Gamma_1 . \,  \Gamma_1 + {\rm O}(\Delta t) \, . 
\monend
In order to compute the coefficient $ \, \Gamma_3 $, we precise now the expansions
of the quantities $ \, A_j \, W + B_j \, Y^* \, $ for $ \, j = 1 $, $\, 2 \, $ and $\, 3 $.
We have

\smallskip  \noindent  $
A \, W + B \, Y^* = 
A \, W + B \, \big( \Phi(W) + \big( \Sigma - {{1}\over{2}} \, {\rm I} \big)   
\, \big( \Delta t \, \Psi_1 +  \Delta t^2 \, \Psi_2 \big) + {\rm O}(\Delta t^3) $

\smallskip  \qquad  \qquad \quad  $ \,\,\,\,
=  A \, W + B \, \Phi + \Delta t \, B \, \big( \Sigma - {{1}\over{2}} \, {\rm I}  \big) \, \Psi_1
+ \Delta t^2 \, B \, \big( \Sigma - {{1}\over{2}} \, {\rm I} \big) \, \Psi_2 
+  {\rm O}(\Delta t^3)  $ 

\smallskip  \noindent  $
A_2 \, W + B_2 \, Y^* =  A_2 \, W + B_2  \, \big( \Phi(W) + \big( \Sigma - {{1}\over{2}} \, {\rm I} \big)
\, \Delta t \, \Psi_1  \big)  + {\rm O}(\Delta t^2) $

\smallskip  \qquad  \qquad \qquad  $ \,\,\, = 
 (A^2 + B\, C) \, W + (A\, B + B \, D ) \, \Phi
+  \Delta t \, B_2 \, \big( \Sigma - {{1}\over{2}} \, {\rm I} \big) \, \Psi_1 + {\rm O}(\Delta t^2) $

\smallskip  \qquad  \qquad \qquad  $ \,\,\, =
A \, \Gamma_1  + B \, ( \gamma_1 - \Psi_1  )  +  \Delta t \, B_2 \, \big( \Sigma - {{1}\over{2}} \, {\rm I} \big) \, \Psi_1
+ {\rm O}(\Delta t^2) $

\hfill because  $ \,\, C \, W + D \, \Psi_1 = \gamma_1 - \Psi_1 $   

\smallskip  \qquad  \qquad \qquad  $ \,\,\, =
\dd \Gamma_1 . \Gamma_1 - B \, \Psi_1 +  \Delta t \, B_2 \, \big( \Sigma - {{1}\over{2}} \, {\rm I} \big) \, \Psi_1
+ {\rm O}(\Delta t^2) $

\hfill because  $ \,\, \dd \Gamma_1 . \Gamma_1 = A \, \Gamma_1 + B \, \gamma_1 $,

\smallskip  \noindent  $
A_3 \, W + B_3 \, Y^* = A_3 \, W + B_3  \,  \Phi(W)  +  {\rm O}(\Delta t) $

\smallskip  \qquad  \qquad \qquad  $ \,\,\, =
(A_2 \, A + B_2 \, C) \, W + (A_2 \, B + B_2 \, D) \, \Phi +  {\rm O}(\Delta t) $

\smallskip  \qquad  \qquad \qquad  $ \,\,\, =
A_2 \, \Gamma_1  + B_2 \, (\gamma_1 - \Psi_1) +  {\rm O}(\Delta t) $

\smallskip  \qquad  \qquad \qquad  $ \,\,\, =
(A^2 + B\, C) \,  \Gamma_1   + (A\, B + B \, D ) \,  \gamma_1 -  B_2 \, \Psi_1  +  {\rm O}(\Delta t) $

\smallskip  \qquad  \qquad \qquad  $ \,\,\, =
A \, (A \, \Gamma_1 + B \, \gamma_1 )  + B \, ( C \,  \Gamma_1  + D  \,  \gamma_1 ) 
-  B_2 \, \Psi_1 +  {\rm O}(\Delta t) $

\smallskip  \qquad  \qquad \qquad  $ \,\,\, = 
A \, \dd \Gamma_1 . \Gamma_1   + B \, ( C \,  \Gamma_1  + D  \,  \gamma_1 ) 
-  B_2 \, \Psi_1 +  {\rm O}(\Delta t) $

\hfill because  $ \,\, \dd \Gamma_1 . \Gamma_1 = A \, \Gamma_1 + B \, \gamma_1 $

\smallskip  \qquad  \qquad \qquad  $ \,\,\, = 
A \, \dd \Gamma_1 . \Gamma_1   + B \,  \dd (\gamma_1 - \Psi_1) . \Gamma_1  
-  B_2 \, \Psi_1 +  {\rm O}(\Delta t) $

\hfill because  $ \,\,
C \, \Gamma_1 + D \, \gamma_1  = \dd (\gamma_1  -\Psi_1) . \Gamma_1 $  

\smallskip  \qquad  \qquad \qquad  $ \,\,\, =
\partial^2 \Gamma_1 . \Gamma_1 -  B \,  \dd \Psi_1 . \Gamma_1   -  B_2 \, \Psi_1  +  {\rm O}(\Delta t) $

\hfill because  $ \,\, \partial^2 \Gamma_1 . \Gamma_1  \equiv \dd ( \dd \Gamma_1 . \Gamma_1 ) . \Gamma_1
=A \, \dd \Gamma_1 . \Gamma_1   + B \,  \dd \gamma_1  . \Gamma_1  $. 

\smallskip  \noindent
We deduce now from (\ref{dvt-W-brut-ordre-3}), (\ref{dt2W-ordre-1}), (\ref{dt3W-ordre-0}) and the previous expressions:

\smallskip  \noindent $ 
\partial_t W =  -A \, W - B \, Y^* +  {1\over2} \, \Delta t \, (A_2 \, W + B_2 \, Y^* - \partial_t^2 W )
- {1\over6} \, \Delta t^2 \, (A_3 \, W + B_3 \, Y^* + \partial_t^3 W ) + {\rm O}(\Delta t^3)  $

\smallskip  \qquad  $ \, = 
 - \big( A \, W + B \, \Phi + \Delta t \, B \, \big( \Sigma - {{1}\over{2}} \, {\rm I}  \big) \, \Psi_1
 + \Delta t^2 \, B \, \big( \Sigma - {{1}\over{2}} \, {\rm I} \big) \, \Psi_2 \big)   $

\smallskip  \qquad  $ \qquad  
+  {1\over2} \, \Delta t \, \big(\dd \Gamma_1 . \Gamma_1 - B \, \Psi_1 +  \Delta t \, B_2 \, \big( \Sigma - {{1}\over{2}} \, {\rm I} \big) \, \Psi_1
- \dd  \Gamma_1 .\, \Gamma_1 - \Delta t \, ( \dd  \Gamma_1 .\, \Gamma_2 + \dd  \Gamma_2 .\, \Gamma_1 ) \big) $

\smallskip  \qquad  $ \qquad  
- {1\over6} \, \Delta t^2 \, \big( \partial^2 \Gamma_1 . \Gamma_1 -  B \,  \dd \Psi_1 . \Gamma_1   -  B_2 \, \Psi_1
- \partial^2 \Gamma_1 . \,  \Gamma_1 \big) + {\rm O}(\Delta t^3) $ 

\smallskip  \qquad  $ \, =
-\Gamma_1 -  \Delta t \, B \, \Sigma \, \Psi_1 + \Delta t^2 \, \big( - B \, \Sigma \, \Psi_2 + {1\over2} \, B \, \Psi_2
+  {1\over2} \,  B_2 \, \Sigma \, \Psi_1  -  {1\over4} \,  B_2 \, \Psi_1  - {1\over2} \, ( \dd  \Gamma_1 .\, \Gamma_2 + \dd  \Gamma_2 .\, \Gamma_1 ) $ 

\smallskip  \qquad  $ \qquad  + 
{1\over6} \, B \, \dd \Psi_1  . \Gamma_1 + {1\over6}  \,  B_2 \, \Psi_1  \big) + {\rm O}(\Delta t^3) $ 
\hfill  because  $ \,\, \partial^2 \Gamma_1 . \Gamma_1  = A \, \dd \Gamma_1 . \Gamma_1   + B \,  \dd \gamma_1  . \Gamma_1  $

\smallskip  \qquad  $ \, =
-\Gamma_1 -  \Delta t \, B \, \Sigma \, \Psi_1 + \Delta t^2 \, \big( - B \, \Sigma \, \Psi_2
+ {1\over2} \, B \, ( \Sigma \,  \dd \Psi_1 . \Gamma_1 + \gamma_2 - D \, \Sigma \,  \Psi_1 ) 
+ {1\over2} \, (A \, B + B \, D) \, \Sigma \, \Psi_1 $

\smallskip  \qquad  $ \qquad  - {1\over12} \, B_2 \, \Psi_1  -  {1\over2} \, (A \, \Gamma_2 + B \, \gamma_2 )  -  {1\over2} \, B \, \Sigma \, \dd \Psi_1 . \Gamma_1   
+ {1\over6} \, B \, \dd \Psi_1  . \Gamma_1 \big)  + {\rm O}(\Delta t^3) $

\smallskip \qquad  $ \, =
-\Gamma_1 -  \Delta t \, \Gamma_2  - \Delta t^2 \, \big( B \, \Sigma \, \Psi_2 +  {1\over12} \, B_2 \, \Psi_1 -  {1\over6} \, B \, \dd \Psi_1  . \Gamma_1  \big)
+ {\rm O}(\Delta t^3) $

\smallskip  \noindent
and the relations (\ref{edp-ordre-3}) and  (\ref{Gamma-3}) are established. \hfill $\square$

\smallskip  \noindent 
In \cite{Du09}, we have proposed a third order expansion of nonlinear lattice Boltzmann schemes for two
specific applications: scalar advection diffusion ($ N= 1$) and fluid system ($ N=d+1$)
with mass and momentum conservation. In the following proposition, we show that our previous
result  at third order accuracy
can be reformulated with the compact formulation (\ref{edp-ordre-3}) (\ref{Gamma-3}) of Proposition~6.

\bigskip \noindent {\bf Proposition 8. \quad   Third-order  formal expansion  for thermics and fluids   } 

\noindent
With the notations (\ref{Lambda-2007}), (\ref{sigma-henon}) and (\ref{defaut-conservation}) 
introduced at Proposition~5, we suppose moreover 
\moneqstar 
M_{0 j} = 1 \,, \,\, M_{\alpha j} = v_j^\alpha \,, \qquad 0 \leq j < q \,,\,\, 1 \leq \alpha \leq d \, . 
\monendstar 
The first moments are denoted by $ \, \rho \equiv \sum_j f_j \, $
and $ \, J_\alpha \equiv \sum_j v_j^\alpha \, f_j  $. 
For the scalar advection diffusion ($ N= 1$), the third order equivalent partial differential  equation  
(\ref{edp-ordre-3}) (\ref{Gamma-3}) can be written
\moneq    \label{thermique-ordre-3}
\partial_t \rho + \partial_\alpha J_\alpha^{\rm eq} - \Delta t \, \widetilde{\sigma}_\alpha \, \partial_\alpha \widetilde{\theta}_\alpha 
+ \Delta t^2 \, \Big[   \Lambda_{\beta \ell}^\gamma \, 
\Big( \widetilde{\sigma}_\beta \,  \widetilde{\sigma}_\ell -  {1\over12} \Big) \,  \partial_\beta \partial_\gamma \widetilde{\theta}_\ell
+  \Big( \widetilde{\sigma}_\beta^2 \,  -  {1\over6} \Big) \, \partial_t  (\partial_\beta \widetilde{\theta}_\beta) \Big]  =   {\rm O}(\Delta t^3)  
\monend
with $ \, J_\alpha^{\rm eq} \equiv \Phi_\alpha (\rho) \, $ and 
the defects of conservation $ \,  \widetilde{\theta}_k (W) \, $ defined in (\ref{defaut-conservation}).
For the  fluid system ($ N=d+1$), we have 
\moneq    \label{fluide-ordre-3} \left\{ \begin{array} {l}
\displaystyle \partial_t \rho + \partial_\alpha J_\alpha  
-  {1\over12} \,\Delta t^2 \,  \Lambda_{\beta \ell}^\gamma \,\,    \partial_\beta \partial_\gamma \widetilde{\theta}_\ell
=  {\rm O}(\Delta t^3) \\ \displaystyle 
\, \partial_t J_\alpha + \Lambda_{\alpha k}^\beta \, \partial_\beta  \, m^{\rm eq}_k 
- \Delta t \,   \Lambda_{\alpha k}^\beta \,\, \widetilde{\sigma}_k \,\, \partial_\beta \, \widetilde{\theta}_k \\ \qquad \qquad \displaystyle
+ \, \Delta t^2  \, \Big[ \, \Lambda_{\alpha k}^\beta \,  \Lambda_{k \ell}^\gamma \, 
\Big(\widetilde{\sigma}_k \,  \widetilde{\sigma}_\ell -  {1\over12} \Big) \,  \partial_\beta \partial_\gamma \widetilde{\theta}_\ell
+ \Lambda_{\alpha k}^\beta \, \Big(\widetilde{\sigma}_k^2 \,  -  {1\over6} \Big) \,  \partial_\beta  \partial_t \widetilde{\theta}_k  \Big] 
 =  {\rm O}(\Delta t^3) \, . 
\end{array} \right. \monend 
 The relations (\ref{thermique-ordre-3}) and (\ref{fluide-ordre-3})
 are axactly the ones proposed with numbers (35), (40) and (41)   in the reference \cite{Du09}.   

\monitem Proof of Proposition 8.  

\noindent
We first consider the conservation defect
$\, \widetilde{\theta}_k(W)  \equiv \partial_t m_k^{\rm eq} + \sum_{\ell \beta} \widetilde{\Lambda}_{k \ell}^\beta \, \partial_\beta  m_\ell^{\rm eq} \, $
introduced in (\ref{defaut-conservation}). We define a vector conservation defect $ \, \theta  \, $ of dimension $ \, q - N \, $
just obtained by a shift of  component numbering:
\moneqstar 
\theta_k =  \widetilde{\theta}_{k+N} (W) \,, \quad 0 \leq k < q-N \, . 
\monendstar 
We have for $\,  0 \leq k < q-N $:

\smallskip \noindent $ 
\theta_k =  ( \dd  \Phi . \, \partial_t W  )_{k} +  ( \widetilde{\Lambda} . m^{\rm eq} )_{k+N}   $

\smallskip \noindent  $ \quad \, =  
( \dd  \Phi . \, \partial_t  W )_k  +  ( C \, W + D \, \Phi(W) )_{k}   $

\smallskip \noindent  $ \quad \, =  
( \dd  \Phi . \, ( -\Gamma_1 - \Delta t \, \Gamma_2 ) )_k   +    ( C \, W + D \, \Phi(W) )_{k} +  {\rm O}(\Delta t^2)  $
\hfill  due to the expansion (\ref{edp-ordre-4}) 

\smallskip \noindent $ \quad \, =  
 - \big( \Psi_1 + \Delta t \,  \dd  \Phi . \, \Gamma_2 \big)_{k}   + {\rm O}(\Delta t^2) \,  $
\hfill  due to the relation  (\ref{formules}) 

\smallskip \noindent
and
\moneqstar 
\theta  =  -  \Psi_1   - \Delta t \,\,   \dd  \Phi . \, \Gamma_2 + {\rm O}(\Delta t^2) \,  .
\monendstar  
 In particular, $ \, \partial_t \theta = - \partial_t \Psi_1  + {\rm O}(\Delta t) =  - \dd \Psi_1 .\, (-\Gamma_1 ) + {\rm O}(\Delta t) $
and
\moneqstar 
\partial_t \theta  =  \dd \Psi_1 . \,\Gamma_1   + {\rm O}(\Delta t) \, . 
\monendstar 
Due to (\ref{Gamma-2}) and (\ref{Gamma-3}), we have the following calculus:

\smallskip \noindent 
$ \, \Gamma_2 + \Delta t \, \Gamma_3 = B \, \Sigma \,  \Psi_1  + \Delta t \, ( B \, \Sigma  \, 
\Psi_2  - {1\over12} \,  B_2 \, \Psi_1  + {{1}\over{6}} \, B \, \dd \Psi_1 . \, \Gamma_1 ) $

\smallskip \noindent  $ \qquad \qquad \quad  =
B \, \Sigma \, (  -\theta +  \Delta t \,\,   \dd  \Phi . \, \Gamma_2 
+ {\rm O}(\Delta t^2)  ) 
+ \Delta t \,  B \, \Sigma \, \big( D \, \Sigma \,  \Psi_1 -   \dd  \Phi . \, \Gamma_2 
- \Sigma \,  \dd \Psi_1 . \,   \Gamma_1  \big) $ 

\smallskip  \noindent  \qquad \qquad  \qquad \qquad  $
- {1\over12} \, \Delta t \, B_2 \, \Psi_1 
+ {{1}\over{6}} \, \Delta t \, B \, \dd \Psi_1 . \, \Gamma_1 $

\smallskip \noindent  $ \qquad \qquad \quad  =
- B \, \Sigma \, \theta + \Delta t \, \big[ 
( B \, \Sigma \,  D \, \Sigma \,  - {1\over12} \, B_2 ) \, \Psi_1
+ B \, (  {{1}\over{6}} - \Sigma^2 ) \, \dd \Psi_1 . \, \Gamma_1 \big]  + {\rm O}(\Delta t^2) $ 

\smallskip \noindent  $ \qquad \qquad \quad  =
- B \, \Sigma \, \theta + \Delta t \, \big[ \, 
- ( B \, \Sigma \,  D \, \Sigma \,  - {1\over12} \,  B_2 ) \, \theta
+ B \, ( \Sigma^2 -  {{1}\over{6}} ) \, \partial_t \theta  \, \big]  + {\rm O}(\Delta t^2) \, . $ 

\smallskip \noindent 
The third-order partial equivalent equations (\ref{edp-ordre-3}) can be written under the form
\moneqstar 
\partial_t W + \widetilde{\Lambda} \, m^{\rm eq} + \Delta t \,  B \, \Sigma \, \theta
+ \Delta t^2 \, \Big[ \, \Big( B \, \Sigma \,  D \, \Sigma \,  - {1\over12} \,  B_2 \Big) \, \theta 
-  B \, \Big( \Sigma^2 -  {{1}\over{6}} \Big) \, \partial_t \theta  \, \Big] =  {\rm O}(\Delta t^3) \, . 
\monendstar  

 \monitem
We explicit now the following cartesian components:

\smallskip \noindent 
$  ( B \, \Sigma \,  D \, \Sigma \, \theta  )_{_{\scriptstyle i}}  = 
- \widetilde{\Lambda}_{i k}^\beta \,\, \widetilde{\sigma}_k \, ( - \widetilde{\Lambda}_{k \ell}^\gamma ) \, \widetilde{\sigma}_\ell \, \, 
\partial_\beta \partial_\gamma \widetilde{\theta}_\ell 
 =  \widetilde{\Lambda}_{i k}^\beta \,\, \widetilde{\sigma}_k \,\, \widetilde{\Lambda}_{k \ell}^\gamma  \, \widetilde{\sigma}_\ell \, \, 
\partial_\beta \partial_\gamma \widetilde{\theta}_\ell  \,, $ 

\smallskip  \noindent 
$  ( B_2 \, \theta  )_{_{\scriptstyle i}}  =   (\widetilde{\Lambda}_{i k}^\beta \, \partial_\beta )\, 
( \widetilde{\Lambda}_{k \ell}^\gamma  \,  \partial_\gamma ) \, \widetilde{\theta}_\ell 
= \widetilde{\Lambda}_{i k}^\beta \,  \widetilde{\Lambda}_{k \ell}^\gamma \, \, \partial_\beta \partial_\gamma \widetilde{\theta}_\ell \,, $ 

\smallskip   \noindent 
$ ( ( B \, \Sigma \,  D \, \Sigma \,  - {1\over12} \,  B_2 ) \, \widetilde{\theta}  )_{_{\scriptstyle i}} =
\widetilde{\Lambda}_{i k}^\beta \,  \widetilde{\Lambda}_{k \ell}^\gamma \, \big(\widetilde{\sigma}_k \,  \widetilde{\sigma}_\ell -  {1\over12} \big) \,
 \partial_\beta \partial_\gamma \widetilde{\theta}_\ell $ 

\smallskip   \noindent    
and $ \,\,  -  ( B \, ( \Sigma^2 -  {{1}\over{6}} ) \, \partial_t \widetilde{\theta} )_{_{\scriptstyle i}} = 
\widetilde{\Lambda}_{i k}^\beta \, \big(\widetilde{\sigma}_k^2 \,  -  {1\over6} \big) \,  \partial_\beta  \partial_t \widetilde{\theta}_k $.  

\smallskip   \noindent
We deduce a new form of the third-order partial equivalent equations for $ \, 0 \leq i < N $: 
\moneq    \label{edp-ordre-3-dcds09-bis} \left\{ \begin{array}{l} 
\partial_t W_i + \widetilde{\Lambda}_{i k}^\alpha \, \partial_\alpha  \, m^{\rm eq}_k 
- \Delta t \,   \widetilde{\Lambda}_{i k}^\beta \, \widetilde{\sigma}_k \, \partial_\beta \, \widetilde{\theta}_k
\\   \vspace {-5mm} \\ 
\qquad  + \, \Delta t^2  \, \big[ \, \widetilde{\Lambda}_{i k}^\beta \,  \widetilde{\Lambda}_{k \ell}^\gamma \, 
\big(\widetilde{\sigma}_k \,  \widetilde{\sigma}_\ell -  {1\over12} \big) \,  \partial_\beta \partial_\gamma \widetilde{\theta}_\ell
+ \widetilde{\Lambda}_{i k}^\beta \, \big(\widetilde{\sigma}_k^2 \,  -  {1\over6} \big) \,  \partial_\beta  \partial_t \widetilde{\theta}_k  \big] 
 =  {\rm O}(\Delta t^3) \, .   
\end {array} \right. \monend

\monitem
We focus on the mass conservation ($ i = 0$). 
We have

\smallskip \noindent $
\widetilde{\Lambda}_{0 \ell}^\alpha  = M_{0 j} \, v_j^\alpha  \, (M^{-1})_{_{\scriptstyle j \ell}} 
=  M_{\alpha j} \, (M^{-1})_{_{\scriptstyle j \ell}} =    \delta_{\alpha \ell}  $. 

\smallskip  \noindent 
Then

\smallskip \noindent $
 \widetilde{\Lambda}_{0 k}^\beta \, \widetilde{\sigma}_k \, \partial_\beta \, \widetilde{\theta}_k = 
 \delta_{\beta k}  \, \widetilde{\sigma}_k \, \partial_\beta \, \widetilde{\theta}_k = 
\widetilde{\sigma}_\beta \, \partial_\beta \, \widetilde{\theta}_\beta  $.  

\smallskip \noindent
In the thermal case, we have only one conservation law and  $ \, \widetilde{\sigma}_\beta \not = 0 $.
Thus we have enlightened  the second order term of the left hand side of the relation (\ref{thermique-ordre-3}). 
In the fluid case, the momentum ({\it id est} moments numbered from $ \, 1 \, $ to $ \, d $) are conserved. In consequence, 
 $ \, \widetilde{\sigma}_\beta  = 0 \, $ and we have no second order term in the mass conservation
(\ref{fluide-ordre-3}). 

\smallskip \noindent
Now  for the first term relative to third order term of (\ref{edp-ordre-3-dcds09-bis}), we have  

\smallskip \noindent $
\widetilde{\Lambda}_{0 k}^\beta \,  \widetilde{\Lambda}_{k \ell}^\gamma \, \big(\widetilde{\sigma}_k \,  \widetilde{\sigma}_\ell -  {1\over12} \big) \,
 \partial_\beta \partial_\gamma \theta_\ell = \delta_{\beta k} \,  \widetilde{\Lambda}_{k \ell}^\gamma \, 
 \big(\widetilde{\sigma}_k \,  \widetilde{\sigma}_\ell -  {1\over12} \big) \,  \partial_\beta \partial_\gamma \theta_\ell
 =  \widetilde{\Lambda}_{\beta \ell}^\gamma \, 
\big(\widetilde{\sigma}_\beta \,  \widetilde{\sigma}_\ell -  {1\over12} \big) \,  \partial_\beta \partial_\gamma \theta_\ell $.  
  
\smallskip \noindent
In the termal case,  $ \, \widetilde{\sigma}_\beta \not = 0 \, $ and this term is exactly the first of the two third order terms
of the relation (\ref{thermique-ordre-3}). In the fluid case,  $ \, \widetilde{\sigma}_\beta  = 0 \, $ and this term has a more compact
expression. Finally,  
\moneq    \label{ordre-3-premier-terme}
\widetilde{\Lambda}_{0 k}^\beta \,  \widetilde{\Lambda}_{k \ell}^\gamma \, \big(\widetilde{\sigma}_k \,  \widetilde{\sigma}_\ell -  {1\over12} \big) \,
 \partial_\beta \partial_\gamma \theta_\ell = 
\left\{ \begin{array} {l}  \widetilde{\Lambda}_{\beta \ell}^\gamma \, 
\big(\widetilde{\sigma}_\beta \,  \widetilde{\sigma}_\ell -  {1\over12} \big) \,  \partial_\beta \partial_\gamma \theta_\ell 
\qquad 
\hfill {\rm thermics} \\ 
-  {1\over12} \,  \widetilde{\Lambda}_{\beta \ell}^\gamma \,   \partial_\beta \partial_\gamma \theta_\ell  \hfill {\rm fluid.} \end{array} 
\right.
\monend

\smallskip \noindent
For the second term relative to  third order, we have 

\smallskip \noindent 
$  \widetilde{\Lambda}_{0 k}^\beta \, \big(\widetilde{\sigma}_k^2 \,  -  {1\over6} \big) \,  \partial_\beta  \partial_t \theta_k =
 \big(\widetilde{\sigma}_\beta^2 \,  -  {1\over6} \big) \, \partial_t  (\partial_\beta \theta_\beta) $

 \smallskip \noindent  
 In the thermal case,  $ \, \widetilde{\sigma}_\beta \not = 0 \, $  and  
$  \widetilde{\Lambda}_{0 k}^\beta \, \big(\widetilde{\sigma}_k^2 \,  -  {1\over6} \big) \,  \partial_\beta  \partial_t \theta_k = 
 \big(\widetilde{\sigma}_\beta^2 \,  -  {1\over6} \big) \, \partial_t  (\partial_\beta \theta_\beta) $.
 In the tluid case,  the momentum is conserved and  $ \, \theta_\beta = {\rm O}(\Delta t)  \, $ and  
$  \, \widetilde{\Lambda}_{0 k}^\beta \, \big(\widetilde{\sigma}_k^2 \,  -  {1\over6} \big) \,  \partial_\beta  \partial_t \theta_k =
 {\rm O}(\Delta t)  $. In a synthetic way,
\moneq    \label{ordre-3-second-terme}
\widetilde{\Lambda}_{0 k}^\beta \, \Big(\widetilde{\sigma}_k^2 \,  -  {1\over6} \Big) \,  \partial_\beta  \partial_t \theta_k =  
\left\{ \begin{array} {l}
 \Big( \widetilde{\sigma}_\beta^2 \,  -  {1\over6} \Big) \, \partial_t  (\partial_\beta \theta_\beta) 
\qquad 
\hfill {\rm thermics} \\ 
 {\rm O}(\Delta t)   \hfill {\rm fluid.} \end{array} 
\right.
\monend  
We deduce from  (\ref{ordre-3-premier-terme}) and (\ref{ordre-3-second-terme})
the expansion  (\ref{thermique-ordre-3}) in the thermal case and the first relation of (\ref{fluide-ordre-3}) in the fluid case.

\monitem
For the momentum equation of the fluid case, $ \, 1 \leq i = \alpha  \leq d \, $ and the relation (\ref{edp-ordre-3-dcds09-bis})
can be rewritten as
\moneqstar  \left\{ \begin{array}{l} 
\partial_t J_\alpha + \widetilde{\Lambda}_{\alpha k}^\alpha \,\, \partial_\alpha  \, m^{\rm eq}_k 
- \Delta t \,   \widetilde{\Lambda}_{\alpha k}^\beta \, \widetilde{\sigma}_k \, \partial_\beta \, \widetilde{\theta}_k
\\   \vspace {-5mm} \\ 
\qquad  + \, \Delta t^2  \, \big[ \, \widetilde{\Lambda}_{\alpha k}^\beta \,  \widetilde{\Lambda}_{k \ell}^\gamma \, 
\big(\widetilde{\sigma}_k \,  \widetilde{\sigma}_\ell -  {1\over12} \big) \,  \partial_\beta \partial_\gamma \widetilde{\theta}_\ell
+ \widetilde{\Lambda}_{\alpha k}^\beta \, \big(\widetilde{\sigma}_k^2 \,  -  {1\over6} \big) \,  \partial_\beta  \partial_t \widetilde{\theta}_k  \big] 
 =  {\rm O}(\Delta t^3) \, .   
\end {array} \right. \monendstar 
This relation is exactly the second relation of (\ref{fluide-ordre-3}).
\hfill $\square$

\bigskip \noindent {\bf Proposition 9. \quad   Nonlinear  expansion at fourth order} 

\noindent
When the time step $ \,  \Delta t \, $ has an infinitesimal value, 
the expansion of the microscopic variables takes the form 
\moneq \label{Y-ordre-3} 
 Y = \Phi(W) +  \Delta t \,\,  S^{-1} \, \Psi_1 (W) +  \Delta t^2 \,\,  S^{-1} \, \Psi_2 (W)  +  \Delta t^3 \,\,  S^{-1} \, \Psi_3 (W)    + {\rm O}(\Delta t^4) \, . 
\monend
The vectors $ \,  \Psi_1 (W) \, $ and  $ \,  \Psi_2 (W) \, $ have been evaluated at the  relations (\ref{Phi-1}) and (\ref{Phi-2}).
We have 
\moneq \label{Phi-3} \left\{ \begin {array} {l}  
\Psi_3  (W)  = \dd \Phi(W) .  \Gamma_3 (W)  +  \Sigma \,\,  \dd \Psi_1 (W) . \Gamma_2 (W) + \Sigma \,\, \dd \Psi_2 (W) . \Gamma_1 (W) - D \, \Sigma \,  \Psi_2 (W) 
\\  \qquad \qquad
+ {1\over6} \, D \, \dd \Psi_1 (W) . \Gamma_1 (W) - {1\over12} \, D_2 \, \Psi_1 (W)
- {1\over12} \, \partial^2  \Psi_1(W) .  \Gamma_1 (W) \, . 
 \end{array}  \right. \monend
A system of fourth order partial differential equations  is asymptotically satisfied by the lattice Boltzmann scheme: 
\moneq \label{edp-ordre-4-final} 
\partial_t W + \Gamma_1 (W) + \Delta t \, \Gamma_2 (W)  + \Delta t^2 \, \Gamma_3 (W)  + \Delta t^3 \, \Gamma_4 (W)  = {\rm O}(\Delta t^4) \, . 
\monend
The vectors  $ \, \Gamma_1 (W) $,  $ \, \Gamma_2 (W) \, $  and  $ \, \Gamma_3 (W) \, $
have been precised at the relations  (\ref{Gamma-1}),   (\ref{Gamma-2})  and (\ref{Gamma-3}) respectively.
We have now 
\moneq \label{Gamma-4}  \left\{ \begin {array} {l}  \displaystyle 
\Gamma_4 (W)  = B \, \Sigma \, \Psi_3 (W) + {1\over4} \, B_2 \, \Psi_2 (W) + {1\over6} \,B \, D_2 \, \Sigma \, \Psi_1 - {1\over6} \, A \, B \, \Psi_2 (W)
\\ \displaystyle \qquad \qquad
- {1\over6} \, B \, \dd \gamma_1 (W) . \Gamma_2(W) - {1\over6} \, B \, \dd \gamma_2 (W) . \Gamma_1 (W) 
- {1\over6} \, B \, \Sigma \,  \partial^2  \Psi_1(W) .  \Gamma_1 (W) \, . 
 \end{array}  \right. \monend

\monitem Proof of Proposition 9. 

\noindent
We first establish the relations (\ref{Y-ordre-3}) and (\ref{Phi-3}). 
We start by the expansion of the scheme (\ref{m-m-star-ordre-3}) at order 3.
For the second component of nonconserved moments, we have

\smallskip \qquad 
$ \left\{ \begin {array} {l}  
Y +  \Delta t \, \partial_t Y + {1\over2} \, \Delta t^2 \, \partial_t^2 Y   + {1\over6} \, \Delta t^3 \, \partial_t^3 Y + {\rm O}(\Delta t^4) =
Y^*  - \Delta t \, (C \, W + D \, Y^*) 
\\  \qquad \qquad         + {1\over2} \, \Delta t^2 \, (C_2 \, W + D_2 \, Y^*)
- {1\over6} \, \Delta t^3 \, (C_3 \, W + D_3 \, Y^*)   + {\rm O}(\Delta t^4) \, . 
 \end{array}  \right. $

\smallskip \noindent 
The relaxation process can be written $ \, Y - Y^*  \equiv  S \, ( Y - \Phi (W) ) $. We deduce
\moneq    \label{dvt-Y-brut-3} \left\{ \begin {array} {l}  \displaystyle 
S \, ( Y - \Phi (W) )  =
- \Delta t \, ( C \, W + D \, Y^* + \partial_t Y) 
+  {1\over2} \, \Delta t^2 \, (C_2 \, W + D_2 \, Y^* - \partial_t^2 Y ) 
\\ \displaystyle \qquad \qquad \qquad \qquad
-  {1\over6} \, \Delta t^3 \, (C_3 \, W + D_3 \, Y^* + \partial_t^3 Y ) + {\rm O}(\Delta t^4) \, . 
 \end{array}  \right. \monend  
We have to  precise the three first time derivatives of the non-conserved moments. 
The relation (\ref{Y-ordre-3}) can also be written at second order accuracy: 
\moneqstar 
Y = \Phi(W) + \Big( \Sigma + {{1}\over{2}} \, {\rm I} \Big)   \,  \Delta t \, \Psi_1 (W)  +
 \Big( \Sigma + {{1}\over{2}} \, {\rm I} \Big)   \,  \Delta t^2 \, \Psi_2 (W) + {\rm O}(\Delta t^3) \, . 
\monendstar  
We differentiate this previous relation relative to time:   

\smallskip \noindent  
$ \, \partial_t Y = \dd \Phi . \, \partial_t W
+ \big( \Sigma + {{1}\over{2}} \, {\rm I} \big)   \,\,   \Delta t \,\,   \dd \Psi_1 . \, \partial_t W
+ \big( \Sigma + {{1}\over{2}} \, {\rm I} \big)   \,\,   \Delta t^2 \,\,   \dd \Psi_2 . \, \partial_t W
+ {\rm O}(\Delta t^3)  $

\smallskip \noindent  \qquad $ = 
- \dd \Phi . \,  \big(  \Gamma_1 (W) + \Delta t \, \Gamma_2 (W) + \Delta t^2 \, \Gamma_3 (W)  + {\rm O}(\Delta t^3)  \big) $

\smallskip \noindent  \qquad \qquad $  
-  \Delta t \,  \big( \Sigma + {{1}\over{2}} \, {\rm I} \big)   \,  \dd \Psi_1 . \,  \big( \Gamma_1  + \Delta t \, \Gamma_2 +  {\rm O}(\Delta t^2)  \big)
+ \, \Delta t^2 \,  \big( \Sigma + {{1}\over{2}} \, {\rm I} \big)   \,  \dd \Psi_2 . \,  \big(  -\Gamma_1  +  {\rm O}(\Delta t)  \big)
+ {\rm O}(\Delta t^3)  $

and 
\moneq    \label{dtY-ordre-2} \left\{ \begin {array} {l} 
\partial_t Y = - \dd \Phi(W)  . \Gamma_1 (W) - \Delta t \, \big[
\dd \Phi(W)  . \Gamma_2 (W) + \big( \Sigma + {{1}\over{2}} \, {\rm I} \big) \,\,   \dd \Psi_1  (W) . \Gamma_1 (W) \big] 
\\ \qquad \qquad 
- \Delta t^2 \, \big[ \dd \Phi(W)  . \Gamma_3 (W) + \big( \Sigma + {{1}\over{2}} \, {\rm I} \big) \,\,   \dd \Psi_1  (W) . \Gamma_2 (W)
\\ \qquad \qquad \qquad 
  + \big( \Sigma + {{1}\over{2}} \, {\rm I} \big) \,\,   \dd \Psi_2  (W) . \Gamma_1 (W) \big] 
+ {\rm O}(\Delta t^3) \,. 
 \end{array}  \right. \monend  
We need the second time derivative $\, \partial_t^2 Y \, $ and third order time derivative $\, \partial_t^3 Y \, $
in the right hand side of (\ref{dvt-Y-brut-3}). We differentiate the relation (\ref{dtY-ordre-2}) relative to time.
Using the notation (\ref{notation-gamma-j}), we deduce

\smallskip \noindent  
$ \, \partial_t^2 Y = -\dd \gamma_1 . (-\Gamma_1 -  \Delta t \, \Gamma_2 )
- \Delta t \, \big[  -\dd \gamma_2 . \Gamma_1 -  \big( \Sigma + {{1}\over{2}} \, {\rm I} \big) \,\,    \partial^2  \Phi(W) .  \Gamma_1 (W)  \big] 
+ {\rm O}(\Delta t^2) $

\noindent and 
\moneq    \label{dt2Y-ordre-1}
\partial_t^2 Y =  \dd \gamma_1 . \, \Gamma_1
+ \, \Delta t \, \Big[ \dd \gamma_1 . \Gamma_2 + \dd \gamma_2 . \Gamma_1
+  \Big( \Sigma + {{1}\over{2}} \, {\rm I} \Big) \,\,    \partial^2  \Psi_1(W) .  \Gamma_1 (W)  \Big] + {\rm O}(\Delta t^2) \, . 
\monend
Finally, 
\moneq    \label{dt3Y-ordre-0}
\partial_t^3 Y = -  \partial^2  \gamma_1(W) .  \Gamma_1 (W) + {\rm O}(\Delta t) \, . 
\monend
We precise now the expansions of $ \,  C \, W + D \, Y^* $, $ \,  C_2 \, W + D_2 \, Y^* \, $ and $ \,  C_3 \, W + D_3 \, Y^* $.
We have 

\smallskip \noindent  
$ C \, W + D \, Y^* =  C \, W + D  \, \big[ \Phi(W) + \big( \Sigma - {{1}\over{2}} \, {\rm I} \big)   \, 
\Delta t \, \Psi_1  + \big( \Sigma - {{1}\over{2}} \, {\rm I} \big) \, \Delta t^2 \, \Psi_2  \big] 
+ {\rm O}(\Delta t^2) $

\noindent  and  
\moneq    \label{CWDYst-ordre-2}  \left\{ \begin {array} {l} 
C \, W + D \, Y^* =  \gamma_1 (W) -  \Psi_1 (W) 
+ \Delta t \,\, D \, \big( \Sigma - {{1}\over{2}} \, {\rm I} \big) \, \Psi_1 (W) 
\\ \qquad \qquad \qquad \qquad 
+ \, \Delta t^2 \,\, D \, \big( \Sigma - {{1}\over{2}} \, {\rm I} \big) \, \Psi_2 (W)   +  {\rm O}(\Delta t^3)  \, . 
 \end{array}  \right. \monend  
For the second order term,  

\smallskip \noindent  
$ C_2 \, W + D_2 \, Y^* =  C_2 \, W + D_2  \, \big[ \Phi(W) +  \Delta t \, \big( \Sigma - {{1}\over{2}} \, {\rm I} \big) \, \Psi_1 \big] 
+ {\rm O}(\Delta t^2)   $

\smallskip and due to (\ref{C2WD2Yst-ordre-0}),
\moneq    \label{C2WD2Yst-ordre-1}
 C_2 \, W + D_2 \, Y^* = \dd ( \gamma_1 - \Psi_1 ) . \Gamma_1 (W) -  D \, \Psi_1  (W) + \Delta t \,  D_2  \,  \Big( \Sigma - {{1}\over{2}} \, {\rm I} \Big) \, \Psi_1 
+ {\rm O}(\Delta t^2)  \, .
\monend  
We have finally for the third order terms

\smallskip \noindent  
$ C_3 \, W + D_3 \, Y^* =   C_3 \, W + D_3  \,\Phi(W) + {\rm O}(\Delta t) $

\smallskip \noindent  \qquad  \qquad  \qquad $ \,\,\,  =
(C_2 \, A + D_2 \, C ) \, W + (C_2 \, B + D_2 \, D ) \, \Phi(W) + {\rm O}(\Delta t) $

\smallskip \noindent  \qquad  \qquad  \qquad $ \,\,\,  =  
C_2 \, \Gamma_1 + D_2 \, ( \gamma_1 - \Psi_1 ) + {\rm O}(\Delta t) $

\smallskip \noindent  \qquad  \qquad  \qquad $ \,\,\,  =  
(C \, A + D \, C ) \, \Gamma_1 +  (C \, B + D^2 ) \, \gamma_1 -   D_2 \, \Psi_1  + {\rm O}(\Delta t) $

\smallskip \noindent  \qquad  \qquad  \qquad $ \,\,\,  =  
C \, \dd \gamma_1 . \, \Gamma_1 + D \, \dd (\gamma_1 - \Psi_1) . \Gamma_1 -   D_2 \, \Psi_1  + {\rm O}(\Delta t) $

\smallskip \noindent  \qquad  \qquad  \qquad $ \,\,\,  =  
C \, \dd \gamma_1 . \, \Gamma_1 + D \, \dd \gamma_1 . \Gamma_1   - D \,  \dd \Psi_1  . \Gamma_1  -   D_2 \, \Psi_1  + {\rm O}(\Delta t) $

\noindent and
\moneq    \label{C3WD3Yst-ordre-0}
C_3 \, W + D_3 \, Y^* = \partial^2 (\gamma_1 - \Psi_1) . \Gamma_1  - D \,   \dd \Psi_1  . \Gamma_1  -   D_2 \, \Psi_1  + {\rm O}(\Delta t)
\monend 
Then due to
(\ref{dvt-Y-brut-3}),  (\ref{dtY-ordre-2}), (\ref{dt2Y-ordre-1}), (\ref{dt3Y-ordre-0}),
(\ref{CWDYst-ordre-2}), (\ref{C2WD2Yst-ordre-1}) and (\ref{C3WD3Yst-ordre-0}), we have

\smallskip \noindent  $
S \, ( Y - \Phi (W) )  =
- \Delta t \, ( C \, W + D \, Y^* + \partial_t Y) 
+  {1\over2} \, \Delta t^2 \, (C_2 \, W + D_2 \, Y^* - \partial_t^2 Y ) $

\smallskip \noindent  \qquad  \qquad \qquad  \qquad $ 
-  {1\over6} \, \Delta t^3 \, (C_3 \, W + D_3 \, Y^* + \partial_t^3 Y ) + {\rm O}(\Delta t^4) $ 

\smallskip \noindent  \quad $ =
- \Delta t \, \big( \gamma_1 -  \Psi_1 + \Delta t \,\, D \, \big( \Sigma - {{1}\over{2}} \, {\rm I} \big) \, \Psi_1 
+ \, \Delta t^2 \,\, D \, \big( \Sigma - {{1}\over{2}} \, {\rm I} \big) \, \Psi_2 \big)   $

\smallskip \noindent  \qquad $ 
-\Delta t \, \big[ - \dd \Phi  . \Gamma_1  - \Delta t \, \big(
\dd \Phi  . \Gamma_2  + \big( \Sigma + {{1}\over{2}} \, {\rm I} \big) \,\,   \dd \Psi_1   . \Gamma_1  \big)  $

\smallskip \noindent  \qquad \qquad $ 
- \Delta t^2 \, \big( \dd \Phi  . \Gamma_3  + \big( \Sigma + {{1}\over{2}} \, {\rm I} \big) \,\,   \dd \Psi_1   . \Gamma_2
  + \big( \Sigma + {{1}\over{2}} \, {\rm I} \big) \,\,   \dd \Psi_2   . \Gamma_1  \big) \big] $

\smallskip \noindent  \qquad $ 
+  {1\over2} \, \Delta t^2 \, \big(
\dd ( \gamma_1 - \Psi_1 ) . \Gamma_1  -  D \, \Psi_1  + \Delta t \,  D_2  \,  \big( \Sigma - {{1}\over{2}} \, {\rm I} \big) \, \Psi_1 \big) $

\smallskip \noindent  \qquad $ 
-  {1\over2} \, \Delta t^2 \, \big( \dd \gamma_1 . \, \Gamma_1
+ \, \Delta t \, \big[ \dd \gamma_1 . \Gamma_2 + \dd \gamma_2 . \Gamma_1
  +  \big( \Sigma + {{1}\over{2}} \, {\rm I} \big) \,\,    \partial^2  \Psi_1 .  \Gamma_1  \big] \big) $

\smallskip \noindent  \qquad $ 
-  {1\over6} \, \Delta t^3 \, \big( \partial^2 (\gamma_1 - \Psi_1) . \Gamma_1  - D \, \dd \Psi_1 . \Gamma_1  -   D_2 \, \Psi_1 \big)
-  {1\over6} \, \Delta t^3 \, \big( -  \partial^2  \gamma_1 .  \Gamma_1 \big) +  {\rm O}(\Delta t^4)  $ 

\newpage 
\smallskip \noindent  \quad $ = 
\Delta t \,\Psi_1  + \Delta t^2 \, \Psi_2 + \Delta t^3 \, \big[
- D \, (\Sigma - {1\over2} \, {\rm I} ) \, \Psi_2  + \dd \Phi . \Gamma_3 +  (\Sigma + {1\over2} \, {\rm I} ) \,  \dd \Psi_1 . \Gamma_2
+  (\Sigma + {1\over2} \, {\rm I} ) \,  \dd \Psi_2 . \Gamma_1 $ 

\smallskip \noindent  \qquad   $ 
+ {1\over2} \, D_2 \,  (\Sigma - {1\over2} \, {\rm I} ) \, \Psi_1
-  {1\over2} \, \dd \gamma_1 . \Gamma_2 -  {1\over2} \, \dd \gamma_2 . \Gamma_1  -  {1\over2} \,   (\Sigma + {1\over2} \, {\rm I} ) \,  \partial^2 \Psi_1 . \Gamma_1
+ {1\over6} \, \partial^2 \Psi_1 . \Gamma_1 + {1\over6} \, D \, \dd \Psi_1 . \Gamma_1 $

\smallskip \noindent  \qquad   $
+  {1\over6} \, D_2 \, \Psi_1  \big] +  {\rm O}(\Delta t^4) $  

\smallskip \noindent  \quad $ =
\Delta t \,\Psi_1  + \Delta t^2 \, \Psi_2 + \Delta t^3 \,  \big[ - D \, \Sigma \, \Psi_2
+ {1\over2} \, D \, ( \Sigma \, \dd \Psi_1 . \Gamma_1 + \gamma_2 - D \, \Sigma \, \Psi_1 ) + \dd \Phi . \Gamma_3
+  \Sigma \, \dd \Psi_1 . \Gamma_2 $

\smallskip \noindent  \qquad   $
+ {1\over2} \,  \dd \Psi_1 . \Gamma_2 +  \Sigma \, \dd \Psi_2 . \Gamma_1
+  {1\over2} \,  ( \Sigma \,\, \partial^2 \Psi_1 . \Gamma_1 + \dd \gamma_2 . \Gamma_1 - D \, \Sigma \, \dd \Psi_1 . \Gamma_1 )
+ {1\over2} \, (C \, B + D^2 ) \, \Sigma \, \Psi_1  $ 

\smallskip \noindent  \qquad   $
- {1\over4} \, D_2 \, \Psi_1 -  {1\over2} \, \dd \gamma_1 . \Gamma_2 -  {1\over2} \, \dd \gamma_2 . \Gamma_1 
-  {1\over2} \, \Sigma \, \, \partial^2 \Psi_1 . \Gamma_1  -  {1\over12} \, \partial^2 \Psi_1 . \Gamma_1
+ {1\over6} \, D \, \dd \Psi_1 . \Gamma_1 +  {1\over6} \, D_2 \, \Psi_1  \big] $

\smallskip \noindent  \hfill   $
+  {\rm O}(\Delta t^4) $

\smallskip \noindent  \quad $ =
\Delta t \,\Psi_1  + \Delta t^2 \, \Psi_2 + \Delta t^3 \,  \big[ \dd \Phi . \Gamma_3 - D \, \Sigma \, \Psi_2 
+ {1\over2} \, D \, \gamma_2 +  \Sigma \, \dd \Psi_1 . \Gamma_2 + {1\over2} \,  \dd \Psi_1 . \Gamma_2 +  \Sigma \,\, \dd \Psi_2 . \Gamma_1  $ 

\smallskip \noindent  \qquad   $
+ \, {1\over2} \, C \, B \, \Sigma \,\, \Psi_1 - {1\over12} \, D_2 \, \Psi_1  - {1\over2} \, \dd \gamma_1 . \Gamma_2
  -  {1\over12} \, \partial^2 \Psi_1 . \Gamma_1 + {1\over6} \, D \, \dd \Psi_1 . \Gamma_1 \big] +  {\rm O}(\Delta t^4)  
$

\smallskip \noindent  \quad $ =
\Delta t \,\Psi_1  + \Delta t^2 \, \Psi_2 + \Delta t^3 \,  \big[ \dd \Phi . \Gamma_3 - D \, \Sigma \, \Psi_2 
+ {1\over2} \, D \, \gamma_2 +  \Sigma \, \dd \Psi_1 . \Gamma_2 + {1\over2} \,  ( -C \, \Gamma_2 - D \, \gamma_2 + \dd \gamma_1 . \Gamma_2 ) $ 

\smallskip \noindent  \qquad   $
+  \Sigma \, \dd \Psi_2 . \Gamma_1 + \, {1\over2} \, C \, B \, \Sigma \,\, \Psi_1 - {1\over12} \, D_2 \, \Psi_1  - {1\over2} \, \dd \gamma_1 . \Gamma_2
  -  {1\over12} \, \partial^2 \Psi_1 . \Gamma_1 + {1\over6} \, D \, \dd \Psi_1 . \Gamma_1 \big] +  {\rm O}(\Delta t^4)   $
  
\smallskip \noindent  \quad $ =
\Delta t \,\Psi_1  + \Delta t^2 \, \Psi_2 + \Delta t^3 \,  \big[ \dd \Phi . \Gamma_3 - D \, \Sigma \, \Psi_2 
+  \Sigma \, \dd \Psi_1 . \Gamma_2 +  \Sigma \, \dd \Psi_2 . \Gamma_1  - {1\over12} \, D_2 \, \Psi_1  
-  {1\over12} \, \partial^2 \Psi_1 . \Gamma_1 $

\smallskip \noindent  \qquad   $
+ {1\over6} \, D \, \dd \Psi_1 . \Gamma_1 \big] +  {\rm O}(\Delta t^4) \, . $ 

\smallskip \noindent  
The relations (\ref{Y-ordre-3}) and (\ref{Phi-3}) are established.

\smallskip \monitem Fourth order partial differential equations

\noindent
We consider now the Taylor expansion of the relation (\ref{exponentielle-Lambda}) at fourth order: 
\moneqstar   \begin {array}{l}  
m + \Delta t \, \partial_t m + {1\over2} \, \Delta t^2 \, \partial_t^2 m + {1\over6} \, \Delta t^3 \, \partial_t^3 m
+ {1\over24} \, \Delta t^4 \, \partial_t^4 m  + {\rm O}(\Delta t^5) 
 \\  \qquad   \qquad  =  m^* - \Delta t \, \Lambda \, m^* + {1\over2} \, \Delta t^2 \, \Lambda^2 \,  m^*   
-{1\over6} \, \Delta t^3 \, \Lambda^3 \,  m^* + {1\over24} \, \Delta t^4 \, \Lambda^4 \,  m^* +  {\rm O}(\Delta t^5)  \, . 
\end{array}  \monendstar  
For the first components (conserved moments):
\moneqstar   \begin {array}{l}
W +  \Delta t \, \partial_t W + {1\over2} \, \Delta t^2 \, \partial_t^2 W 
+ {1\over6} \, \Delta t^3 \, \partial_t^3 W + {1\over24} \, \Delta t^4 \, \partial_t^4 W + {\rm O}(\Delta t^5)
=  W - \Delta t \, (A \, W + B \, Y^*) 
 \\  \qquad      + {1\over2} \, \Delta t^2 \, (A_2 \, W + B_2 \, Y^*)
-  {1\over6} \, \Delta t^3 \, (A_3 \, W + B_3 \, Y^*)  +  {1\over24} \, \Delta t^4 \, (A_4 \, W + B_4 \, Y^*)  + {\rm O}(\Delta t^5) \,. 
\end{array}  \monendstar   
We simplify the constant term $ \, W \, $ and divide by $ \, \Delta t $. We deduce
\moneq \label{dvt-W-brut-ordre-4}   \left \{ \begin {array}{l}  
\partial_t W =  - (A \, W + B \, Y^*) +  {1\over2} \, \Delta t \, (A_2 \, W + B_2 \, Y^* - \partial_t^2 W )
\\  \qquad  - {1\over6} \, \Delta t^2 \, (A_3 \, W + B_3 \, Y^* + \partial_t^3 W )
+ {1\over24} \, \Delta t^3 \, (A_4 \, W + B_4 \, Y^* - \partial_t^4 W ) + {\rm O}(\Delta t^4)  \,. 
\end{array} \right. \monend
We explicit the partial derivatives $ \, \partial_t^2 W  $, $ \, \partial_t^3 W \, $  and $ \, \partial_t^4 W \, $  at orders two,
one and zero respectively. We start with the relation 
\moneqstar
\partial_t W + \Gamma_1 + \Delta t \, \Gamma_2   + \Delta t^2 \, \Gamma_3  = {\rm O}(\Delta t^3)  
\monendstar 
where $ \, \Gamma_1 $, $ \, \Gamma_2 \, $ and $ \, \Gamma_3 \, $ have been evaluated previously. 
We have

\smallskip \noindent  $
\partial_t^2 W = \partial_t \big( -\Gamma_1 - \Delta t \, \Gamma_2  - \Delta t^2 \, \Gamma_3 + {\rm O}(\Delta t^3)  \big)   $ 

\smallskip \noindent  \qquad  $  \,\, = 
\dd \big( -\Gamma_1 - \Delta t \, \Gamma_2  - \Delta t^2 \, \Gamma_3 + {\rm O}(\Delta t^3) \big) . (\partial_t W) $ 

\smallskip \noindent  \qquad  $  \,\, =
\dd \big( -\Gamma_1 - \Delta t \, \Gamma_2 - \Delta t^2 \, \Gamma_3 + {\rm O}(\Delta t^2) \big) \, . \, \big( -\Gamma_1 - \Delta t \, \Gamma_2 - \Delta t^2 \, \Gamma_3 \big) 
+ {\rm O}(\Delta t^3)  $

\smallskip \noindent
and 
\moneq    \label{dt2W-ordre-2}
\partial_t^2 W = 
\dd  \Gamma_1 .\, \Gamma_1 + \Delta t \, ( \dd  \Gamma_1 .\, \Gamma_2 + \dd  \Gamma_2 .\, \Gamma_1 )
+ \Delta t^2 \, ( \dd  \Gamma_1 .\, \Gamma_3 + \dd  \Gamma_2 .\, \Gamma_2 + \dd  \Gamma_3 .\, \Gamma_1 ) + {\rm O}(\Delta t^2) \, . 
\monend
Then 

\smallskip \noindent 
$\, \partial_t^3 W = \partial_t (  \partial_t^2 W )  $

\smallskip  \qquad 
$  \,\,\, =  \partial_t \big( A \,  \Gamma_1 + B \,  \dd \Phi . \Gamma_1  
+ \Delta t \, \partial_t ( \dd  \Gamma_2 .\, \Gamma_1 + \dd  \Gamma_1 .\, \Gamma_2 ) \big) + {\rm O}(\Delta t^2) $

\smallskip  \qquad 
$ \,\,\, =  ( A \, \dd \Gamma_1 + B \, \dd  ( \dd \Phi . \Gamma_1 ) ) . \, \big( - \Gamma_1 
- \Delta t \, \Gamma_2 + {\rm O}(\Delta t^2) \big) $

\qquad \qquad 
$ +    \Delta t \, \partial_t \, ( B \, \Sigma \, \dd \Psi_1 . \, \Gamma_1 + A \, \Gamma_2  + B \, \dd \Phi . \, \Gamma_2 ) 
 + {\rm O}(\Delta t^2) $ 

\smallskip  \qquad 
$  \,\,\, =  \partial_t \big( \dd  \Gamma_1 .\, \Gamma_1 + \Delta t \, ( \dd  \Gamma_1 .\, \Gamma_2 + \dd  \Gamma_2 .\, \Gamma_1 ) \big)
+ {\rm O}(\Delta t^2) $  

\smallskip  \qquad 
$ \,\,\, = \dd \big( \dd  \Gamma_1 .\, \Gamma_1+ \Delta t \, ( \dd  \Gamma_1 .\, \Gamma_2 + \dd  \Gamma_2 .\, \Gamma_1 ) \big) . \, \partial_t W  + {\rm O}(\Delta t^2) $  

\smallskip  \qquad 
$ \,\,\, = \dd \big( \dd  \Gamma_1 .\, \Gamma_1+ \Delta t \, ( \dd  \Gamma_1 .\, \Gamma_2 + \dd  \Gamma_2 .\, \Gamma_1 ) \big) .
\big( -\Gamma_1  - \Delta t \, \Gamma_2 ) + {\rm O}(\Delta t^2) $

\smallskip  \qquad 
$ \,\,\, =  -\dd \big( \dd  \Gamma_1 .\, \Gamma_1 \big) . \Gamma_1 - \Delta t \, \big[ \dd ( A \, \Gamma_2 + B \, \gamma_2
+ B \, \Sigma \, \dd  \Psi_1 .\, \Gamma_1 ) . \, \Gamma_1 + \dd (  A \, \Gamma_1 + B \, \gamma_1 ) . \, \Gamma_2  \big]  + {\rm O}(\Delta t^2) $

\smallskip  \qquad 
$ \,\,\, =  -\partial^2 \Gamma_1 . \,  \Gamma_1  - \Delta t \, \big[ A \, \dd \Gamma_2 . \, \Gamma_1 + B \, \dd \gamma_2  . \, \Gamma_1
+ B \, \Sigma \, \partial^2 \Psi_1 . \,  \Gamma_1  +  A \, \dd \Gamma_1 . \, \Gamma_2 +  B \, \dd \gamma_1  . \, \Gamma_2 \big] 
+ {\rm O}(\Delta t^2) $  

and we have at first order for this third order derivative
\moneq    \label{dt3W-ordre-1} \left \{ \begin {array}{l}  
\partial_t^3 W = -\partial^2 \Gamma_1 . \,  \Gamma_1  - \Delta t \, \big[ A \, (\dd \Gamma_2 . \, \Gamma_1  + \dd \Gamma_1 . \, \Gamma_2 ) 
\\  \qquad \qquad \qquad 
+ B \, (\dd \gamma_2  . \, \Gamma_1 + \dd \gamma_1  . \, \Gamma_2 +  \Sigma \, \partial^2 \Psi_1 . \,  \Gamma_1  )  \big]  + {\rm O}(\Delta t^2) \, . 
\end{array} \right. \monend
Last but not least, we have

\smallskip \noindent 
$  \partial_t^4 W = \partial_t \big(  -\partial^2 \Gamma_1 . \,  \Gamma_1 + {\rm O}(\Delta t) \big) $

\smallskip  \qquad 
$ \,\, = - \partial_t \big( A \, \dd \Gamma_1 . \, \Gamma_1 + B \, \dd \gamma_1 . \, \Gamma_1 )  + {\rm O}(\Delta t) \big) $

\smallskip  \qquad 
$ \,\, = \dd \big(  A \, \dd \Gamma_1 . \, \Gamma_1 + B \, \dd \gamma_1 . \, \Gamma_1 \big) . \, ( \Gamma_1  + {\rm O}(\Delta t) ) \,  $ 

\smallskip \noindent
and finally 
\moneq    \label{dt4W-ordre-0} 
 \partial_t^4 W =  A \, \partial^2 \Gamma_1 . \,  \Gamma_1 + B \,  \partial^2 \gamma_1 . \,  \Gamma_1  + {\rm O}(\Delta t)  \, . 
 \monend
For the determination of the coefficient $ \, \Gamma_4 $, we precise  the expansions
of the quantities $ \, A_j \, W + B_j \, Y^* \, $ for $ \, j = 1 \, $ to $\, j=4 $.
We use the previous evaluations done in the proof of the previous propositions: 

\smallskip  \noindent  $
A \, W + B \, Y^* = 
A \, W + B \, \big[ \Phi(W) + \big( \Sigma - {{1}\over{2}} \, {\rm I} \big)   
\, \big( \Delta t \, \Psi_1 +  \Delta t^2 \, \Psi_2  +  \Delta t^3 \, \Psi_3 \big) \big] + {\rm O}(\Delta t^4) $

\smallskip  \qquad  \qquad \quad  $ \,\,\,\,
=  A \, W + B \, \Phi + \Delta t \, B \, \big( \Sigma - {{1}\over{2}} \, {\rm I}  \big) \, \Psi_1
+ \Delta t^2 \, B \, \big( \Sigma - {{1}\over{2}} \, {\rm I} \big) \, \Psi_2 
+ \Delta t^3 \, B \, \big( \Sigma - {{1}\over{2}} \, {\rm I} \big) \, \Psi_3 
+  {\rm O}(\Delta t^4)  $ 

\smallskip  \noindent  $
A_2 \, W + B_2 \, Y^* =  A_2 \, W + B_2  \, \big[ \Phi(W) + \big( \Sigma - {{1}\over{2}} \, {\rm I} \big)
  \, ( \Delta t \, \Psi_1 + \Delta t^2 \, \Psi_2 )    \big]
  + {\rm O}(\Delta t^3) $

\smallskip  \qquad  \qquad  $ = 
 (A^2 + B\, C) \, W + (A\, B + B \, D ) \, \Phi
+  \Delta t \, B_2 \, \big( \Sigma - {{1}\over{2}} \, {\rm I} \big) \, \Psi_1 
+  \Delta t^2 \, B_2 \, \big( \Sigma - {{1}\over{2}} \, {\rm I} \big) \, \Psi_2 + {\rm O}(\Delta t^3) $

\smallskip  \qquad  \qquad  $ = \dd \Gamma_1 . \Gamma_1 - B \, \Psi_1 +  \Delta t \, B_2 \, \big( \Sigma - {{1}\over{2}} \, {\rm I} \big) \, \Psi_1
+  \Delta t^2 \, B_2 \, \big( \Sigma - {{1}\over{2}} \, {\rm I} \big) \, \Psi_2 + {\rm O}(\Delta t^3) $

\smallskip  \noindent  $
A_3 \, W + B_3 \, Y^* = A_3 \, W + B_3  \, [ \Phi(W) +  \Delta t \, \big( \Sigma - {{1}\over{2}} \, {\rm I} \big)   
\, \Psi_1 ]   +  {\rm O}(\Delta t^2) $

\smallskip  \qquad  \qquad \qquad  $ \,\,\, =
(A_2 \, A + B_2 \, C) \, W + (A_2 \, B + B_2 \, D) \, \Phi +
 \Delta t \, B_3 \,  \big( \Sigma - {{1}\over{2}} \, {\rm I} \big) \, \Psi_1 
+ {\rm O}(\Delta t^2) $

\smallskip  \qquad  \qquad \qquad  $ \,\,\, =
\partial^2 \Gamma_1 . \Gamma_1 -  B \,  \dd \Psi_1 . \Gamma_1   -  B_2 \, \Psi_1
+  \Delta t \, B_3 \,  \big( \Sigma - {{1}\over{2}} \, {\rm I} \big) \, \Psi_1 
+  {\rm O}(\Delta t^2) $

\smallskip  \noindent  $
A_4 \, W + B_4 \, Y^* = A_4 \, W + B_4  \,  \Phi(W) +  {\rm O}(\Delta t) $ 

\smallskip  \qquad  \qquad \qquad  $ \,\,\, =
(A_3 \, A + B_3 \, C ) \, W + ( A_3 \, B + B_3 \, D ) \, \Phi  +  {\rm O}(\Delta t) $ 

\smallskip  \qquad  \qquad \qquad  $ \,\,\, =
A_3 \, \Gamma_1 + B_3 \, ( \gamma_1 - \Psi_1 )   +  {\rm O}(\Delta t) $ 

\smallskip  \qquad  \qquad \qquad  $ \,\,\, = 
(A_2 \, A + B_2 \, C ) \,\Gamma_1 + ( A_2 \, B + B_2 \, D ) \,\gamma_1 -  B_3 \, \Psi_1  +  {\rm O}(\Delta t) $ 

\smallskip  \qquad  \qquad \qquad  $ \,\,\, = 
A_2 \, \dd \Gamma_1 . \, \Gamma_1 + B_2 \,  \dd (\gamma_1 - \Psi_1 ) . \, \Gamma_1 -  B_3 \, \Psi_1  +  {\rm O}(\Delta t) $ 

\smallskip  \qquad  \qquad \qquad  $ \,\,\, =
(A^2 + B \, C ) \, \dd \Gamma_1 . \, \Gamma_1 + (A \, B + B \, D ) \, \dd \gamma_1 . \, \Gamma_1 - B_2 \, \dd \Psi_1  . \, \Gamma_1
-  B_3 \, \Psi_1  +  {\rm O}(\Delta t) $

\smallskip  \qquad  \qquad \qquad  $ \,\,\, =
A \, \partial^2 \Gamma_1 . \, \Gamma_1 + B \, \partial^2 (\gamma_1 - \Psi_1) . \, \Gamma_1
- B_2 \, \dd \Psi_1  . \, \Gamma_1 -  B_3 \, \Psi_1  +  {\rm O}(\Delta t)  $

\smallskip  \qquad  \qquad \qquad  $ \,\,\, =
A \, \partial^2 \Gamma_1 . \, \Gamma_1 + B \, \partial^2 \gamma_1 . \, \Gamma_1 - B \,  \partial^2 \Psi_1 . \, \Gamma_1
- B_2 \, \dd \Psi_1  . \, \Gamma_1 -  B_3 \, \Psi_1  +  {\rm O}(\Delta t) \, . $ 

\smallskip  \noindent
We deduce now from (\ref{dvt-W-brut-ordre-4}), (\ref{dt2W-ordre-2}), (\ref{dt3W-ordre-1}), (\ref{dt4W-ordre-0}) 
and the previous expressions:

\smallskip  \noindent $ 
\partial_t W =  - (A \, W + B \, Y^*) +  {1\over2} \, \Delta t \, (A_2 \, W + B_2 \, Y^* - \partial_t^2 W )
- {1\over6} \, \Delta t^2 \, (A_3 \, W + B_3 \, Y^* + \partial_t^3 W ) $

\smallskip  \qquad  \qquad  $ 
+ {1\over24} \, \Delta t^3 \, (A_4 \, W + B_4 \, Y^* - \partial_t^4 W )
+ {\rm O}(\Delta t^4)  $

\smallskip  \qquad  $ \, =
 - \big(  A \, W + B \, \Phi + \Delta t \, B \, \big( \Sigma - {{1}\over{2}} \, {\rm I}  \big) \, \Psi_1
+ \Delta t^2 \, B \, \big( \Sigma - {{1}\over{2}} \, {\rm I} \big) \, \Psi_2 
+ \Delta t^3 \, B \, \big( \Sigma - {{1}\over{2}} \, {\rm I} \big) \, \Psi_3 \big)   $

\smallskip  \qquad  $ \qquad  
+  {1\over2} \, \Delta t \, \big[
\dd \Gamma_1 . \Gamma_1 - B \, \Psi_1 +  \Delta t \, B_2 \, \big( \Sigma - {{1}\over{2}} \, {\rm I} \big) \, \Psi_1
+  \Delta t^2 \, B_2 \, \big( \Sigma - {{1}\over{2}} \, {\rm I} \big) \, \Psi_2 $ 

\smallskip  \qquad  $ \qquad  \quad  
- \big(  \dd  \Gamma_1 .\, \Gamma_1 + \Delta t \, ( \dd  \Gamma_1 .\, \Gamma_2 + \dd  \Gamma_2 .\, \Gamma_1 )
+ \Delta t^2 \, ( \dd  \Gamma_1 .\, \Gamma_3 + \dd  \Gamma_2 .\, \Gamma_2 + \dd  \Gamma_3 .\, \Gamma_1 ) \big) \big] $

\smallskip  \qquad  $ \qquad  
- {1\over6} \, \Delta t^2 \, \big[
\partial^2 \Gamma_1 . \Gamma_1 -  B \,  \dd \Psi_1 . \Gamma_1   -  B_2 \, \Psi_1
+  \Delta t \, B_3 \,  \big( \Sigma - {{1}\over{2}} \, {\rm I} \big) \, \Psi_1 $

\smallskip  \qquad  $ \qquad  \quad
 -\partial^2 \Gamma_1 . \,  \Gamma_1  - \Delta t \, \big( A \, (\dd \Gamma_2 . \, \Gamma_1  + \dd \Gamma_1 . \, \Gamma_2 ) 
+ B \, (\dd \gamma_2  . \, \Gamma_1 + \dd \gamma_1  . \, \Gamma_2 +  \Sigma \, \partial^2 \Psi_1 . \,  \Gamma_1  )   \big) \big] $

\smallskip  \qquad  $ \qquad  
+ {1\over24} \, \Delta t^3 \, \big[
A \, \partial^2 \Gamma_1 . \, \Gamma_1 + B \, \partial^2 \gamma_1 . \, \Gamma_1 - B \,  \partial^2 \Psi_1 . \, \Gamma_1
- B_2 \, \dd \Psi_1  . \, \Gamma_1 -  B_3 \, \Psi_1 $

\smallskip  \qquad  $ \qquad  \quad
- \big(  A \, \partial^2 \Gamma_1 . \,  \Gamma_1 + B \,  \partial^2 \gamma_1 . \,  \Gamma_1 \big) \big]   + {\rm O}(\Delta t^4) $

\smallskip  \qquad  $ \, =
-\Gamma_1 -  \Delta t \, \Gamma_2 - \Delta t^2 \, \Gamma_3 +  \Delta t^3 \, \big[  - B \, \big( \Sigma - {{1}\over{2}} \, {\rm I} \big) \, \Psi_3
+  {1\over2} \,  B_2 \, \big( \Sigma - {{1}\over{2}} \, {\rm I} \big) \, \Psi_2 $

\smallskip  \qquad  $ \qquad 
-  {1\over2} \, ( \dd  \Gamma_1 .\, \Gamma_3 + \dd  \Gamma_2 .\, \Gamma_2 + \dd  \Gamma_3 .\, \Gamma_1 )
+  {1\over6} \,  \big( -B_3 \, \big( \Sigma - {{1}\over{2}} \, {\rm I} \big) \, \Psi_1
+ A \, (\dd \Gamma_2 . \, \Gamma_1  + \dd \Gamma_1 . \, \Gamma_2 ) $

\smallskip  \qquad  $ \qquad 
+ B \, (\dd \gamma_2  . \, \Gamma_1 + \dd \gamma_1  . \, \Gamma_2 +  \Sigma \, \partial^2 \Psi_1 . \,  \Gamma_1  )   \big)
+  {1\over24} \,  \big( A \, \partial^2 \Gamma_1 . \, \Gamma_1 + B \, \partial^2 \gamma_1 . \, \Gamma_1 $

\smallskip  \qquad  $ \qquad
- B \,  \partial^2 \Psi_1 . \, \Gamma_1 - B_2 \, \dd \Psi_1  . \, \Gamma_1 -  B_3 \, \Psi_1
- A \, \partial^2 \Gamma_1 . \,  \Gamma_1 - B \,  \partial^2 \gamma_1 . \,  \Gamma_1 \big) \big] +  {\rm O}(\Delta t^4) $ 

\smallskip  \qquad  $ \, =
-\Gamma_1 -  \Delta t \, \Gamma_2 - \Delta t^2 \, \Gamma_3  -  \Delta t^3 \, \big[ 
 B \, \Sigma  \, \Psi_3 -  {1\over2} \,  B \, \big( \gamma_3 + \Sigma \, \dd \Psi_1. \Gamma_2 + \Sigma \, \dd \Psi_2 . \Gamma_1 - D \, \Sigma \, \Psi_2 $

\smallskip  \qquad  $ \qquad
+ {1\over6} \,  D \, \dd \Psi_1 . \Gamma_1 - {1\over12} \, D_2 \, \Psi_1 - {1\over12} \, \partial^2  \Psi_1 .  \Gamma_1  \big) 
- {1\over2} \, (A \, B + B \, D) \, \Sigma \, \Psi_2  + {1\over4} \,  B_2 \, \Psi_2  $

\smallskip  \qquad  $ \qquad
 + {1\over2} \, \big( A \, \Gamma_3 + B \, \gamma_3 + B \, \Sigma \, \dd \Psi_1 . \Gamma_1 +  \dd \Gamma_3 . \, \Gamma_1 \big) 
+ {1\over6} \, B_3 \, \Sigma \, \Psi_1  -  {1\over12} \, B_3  \, \Psi_1 - {1\over6} \, A \, B \, \Sigma \,  \dd \Psi_1 .  \Gamma_1 $

\smallskip  \qquad  $ \qquad
  - {1\over6} \, A \, ( A \, \Gamma_2 + B \, \gamma_2)  
  - {1\over6} \, B \, \big( \dd \gamma_2  . \, \Gamma_1 + \dd \gamma_1  . \, \Gamma_2 +  \Sigma \, \partial^2 \Psi_1 . \,  \Gamma_1    \big)
+  {1\over24} \, B \,  \partial^2 \Psi_1 . \,  \Gamma_1 +  {1\over24} \, B_2 \, \dd \Psi_1  . \, \Gamma_1 $

\smallskip  \qquad  $ \qquad 
+  {1\over24} \,B_3 \, \Psi_1   \big] +   {\rm O}(\Delta t^4) $

\smallskip  \qquad  $ \, =
-\Gamma_1 -  \Delta t \, \Gamma_2 - \Delta t^2 \, \Gamma_3  -  \Delta t^3 \, \big[ 
 B \, \Sigma  \, \Psi_3 - {1\over2} \, B \, \Sigma \, \dd \Psi_2 . \Gamma_1    - {1\over12} \, B \, D \,  \dd \Psi_1 . \Gamma_1
+  {1\over24} \, B \, D_2 \, \Psi_1 $ 

\smallskip  \qquad  $ \qquad
+  {1\over24} \, B \, \partial^2  \Psi_1 .  \Gamma_1  - {1\over2} \, A \, B \, \Sigma \, \Psi_2  + {1\over4} \,  B_2 \, \Psi_2 
+ {1\over2} \, A \, \big( B \, \Sigma \, \Psi_2 -   {1\over6} \, B \,  \dd \Psi_1 . \Gamma_1 +  {1\over12} \, B_2 \, \Psi_1 \big) $

\smallskip  \qquad  $ \qquad
+  {1\over2} \, \big( B \, \Sigma \, \dd \Psi_2 . \Gamma_1 -   {1\over6} \, B \,  \partial^2 \Psi_1 . \Gamma_1  +  {1\over12} \, B_2 \, \dd \Psi_1  . \Gamma_1 \big)
+ {1\over6} \, B_3 \, \Sigma \, \Psi_1  - {1\over24} \,B_3 \, \Psi_1 - {1\over6} \, A \, B \, \Sigma \,  \dd \Psi_1 .  \Gamma_1  $

\smallskip  \qquad  $ \qquad
-  {1\over6} \, A^2 \, \Gamma_2 - {1\over6} \, A \, B \, \gamma_2 
 - {1\over6} \, B \, \big( \dd \gamma_2  . \, \Gamma_1 + \dd \gamma_1  . \, \Gamma_2 +  \Sigma \, \partial^2 \Psi_1 . \,  \Gamma_1    \big)
+  {1\over24} \, B \,  \partial^2 \Psi_1  . \,  \Gamma_1  $

\smallskip  \qquad  $ \qquad
+  {1\over24} \, B_2 \, \dd \Psi_1  . \, \Gamma_1  \big] +   {\rm O}(\Delta t^4) $

\smallskip  \qquad  $ \, =
-\Gamma_1 -  \Delta t \, \Gamma_2 - \Delta t^2 \, \Gamma_3  -  \Delta t^3 \, \big[ 
 B \, \Sigma  \, \Psi_3 -  {1\over12} \, B_2 \, \dd \Psi_1  . \, \Gamma_1 +  {1\over24} \, B \, D_2 \, \Psi_1 +  {1\over4} \,  B_2 \, \Psi_2 
+  {1\over24} \, A \, B_2 \, \Psi_1  $

\smallskip  \qquad  $ \qquad
+ {1\over6} \, B_3 \, \Sigma \, \Psi_1  - {1\over24} \,B_3 \, \Psi_1  -  {1\over6} \, A \, B \, \Sigma \, \dd \Psi_1 .  \Gamma_1 -  {1\over6} \, A^2 \, B  \, \Sigma \, \Psi_1
 -  {1\over6} \, A \, B \, \gamma_2  $

\smallskip  \qquad  $ \qquad
 - {1\over6} \, B \, \big( \dd \gamma_2  . \, \Gamma_1 + \dd \gamma_1  . \, \Gamma_2 +  \Sigma \, \partial^2 \Psi_1 . \,  \Gamma_1    \big) 
 +  {1\over12} \, B_2 \, \dd \Psi_1  . \, \Gamma_1  \big] +   {\rm O}(\Delta t^4) $

\smallskip  \qquad  $ \, =
-\Gamma_1 -  \Delta t \, \Gamma_2 - \Delta t^2 \, \Gamma_3  -  \Delta t^3 \, \big[ 
 B \, \Sigma  \, \Psi_3  +  {1\over4} \,  B_2 \, \Psi_2 + {1\over6} \, B_3 \, \Sigma \, \Psi_1 
 -  {1\over6} \, A \, B \, \Sigma \, \dd \Psi_1 .  \Gamma_1 $

 \smallskip  \qquad  $ \qquad
 -  {1\over6} \, A^2 \, B  \, \Sigma \, \Psi_1  -  {1\over6} \, A \, B \, \big( \Psi_2 - \Sigma \, \dd \Psi_1  . \, \Gamma_1 + D \, \Sigma \, \Psi_1 \big) 
 - {1\over6} \, B \, \big( \dd \gamma_2  . \, \Gamma_1 + \dd \gamma_1  . \, \Gamma_2 +  \Sigma \, \partial^2 \Psi_1 . \,  \Gamma_1    \big) \big] $

\smallskip  \qquad  $ \qquad 
+  {\rm O}(\Delta t^4) $
\hfill because $ \, A \, B_2 + B \, D_2 = B_3 $ 

\smallskip  \qquad  $ \, =
-\Gamma_1 -  \Delta t \, \Gamma_2 - \Delta t^2 \, \Gamma_3  -  \Delta t^3 \, \big[ 
 B \, \Sigma  \, \Psi_3  +  {1\over4} \,  B_2 \, \Psi_2   + {1\over6} \, B \, D_2 \, \Sigma \, \Psi_1 -   {1\over6} \, A \, B \,  \Psi_2  $

\smallskip  \qquad  $ \qquad
 - {1\over6} \, B \, \big( \dd \gamma_2  . \, \Gamma_1 + \dd \gamma_1  . \, \Gamma_2 +  \Sigma \, \partial^2 \Psi_1 . \,  \Gamma_1    \big) \big] 
+  {\rm O}(\Delta t^4)  $

\smallskip \noindent 
 because $ \, B_3 - A \, (A \, B +  B \, D) = B \, D_2 $. 
The relations (\ref{edp-ordre-4-final}) and  (\ref{Gamma-4}) are established
and the proposition is proven. \hfill $\square$

\newpage 
\bigskip \bigskip     \noindent {\bf \large  7) \quad  Revisiting the ``Berlin algorithm''  in the linear case }    

In \cite{ADGL14}, we have presented the Berlin explicit algorithm  in order to determine, with the help of formal calculus, 
the equivalent partial differential equations of a linear lattice Boltzmann scheme.  
We do not enter here into the details of our previous contribution because
the linearized version of the present algorithm is much more simple. 

\bigskip \noindent {\bf Proposition 10. \quad   Linearized general  expansion at fourth order } 

We suppose that the equilibrium value of the nonconserved moments 
is a linear function:
\moneqstar 
 \Phi(W) =  E \,\, W  \,  
\monendstar 
with a given fixed $ \, (q-N) \times N \, $ rectangle matrix $ \, E $. 
Then the nonequilibium moments can be expanded as 
\moneqstar 
Y = \big[ E + S^{-1} ( \Delta t \, \beta_1 + \Delta t^2 \, \beta_2 + \Delta t^3 \, \beta_3 ) \big] \, W + {\rm O}(\Delta t^4)  \, . 
\monendstar 
The  $ \, (q-N) \times N \, $ matrices $ \, \beta_j \, $  are linear operators of order $ \, j $. 
The equivalent partial equivalent equations are linear and can be written as 
\moneqstar 
\partial_t W + \big[  \alpha_1 + \Delta t \, \alpha_2 + \Delta t^2 \, \alpha_3  + \Delta t^3 \, \alpha_4 \big] \, W + {\rm O}(\Delta t^4)  \, , 
\monendstar 
and $ \, \alpha_j \, $ is a  $ \, N \times N \, $ differential  matrix of order $ \, j $. 
In other terms, $ \, \Phi_j (W) = \beta_j \, W \, $ and $ \, \Gamma_j(W) = \alpha_j \, W \, $ are linear functions 
of the conserved moments. 
The coefficent matrices $ \,  \alpha_j \, $ and $ \,  \beta_j \, $ are computed according to 
a linearized version of (\ref{formules}):
\moneq \label{formules-lineaires} \left \{ \begin {array}{rl}
\alpha_1   & \!\!\! =     A  + B \, E \\
\beta_1  & \!\!\! =   E \, \alpha_1 - C - D \, E \\ 
\alpha_2  & \!\!\! =  B \, \Sigma \, \beta_1  \\
\beta_2  & \!\!\! =   \Sigma \,   \beta_1      \alpha_1  +  E \, \alpha_2  
- D \, \Sigma \,  \beta_1    \\ 
\alpha_3   & \!\!\! = B \, \Sigma  \, \beta_2  + {{1}\over{12}}  B_2 \, \beta_1  
-  {{1}\over{6}} \, B \,   \beta_1 \, \alpha_1  \\
\beta_3  & \!\!\! = \Sigma \,  \beta_1     \alpha_2   +   E \, \alpha_3 -  D \, \Sigma \, \beta_2  
 + \Sigma \,  \beta_2     \alpha_1   +{1\over6} \, D \,  \beta_1     \alpha_1   - {1\over12} \, D_2 \, \beta_1   
- {1\over12} \, \beta_1 \, \alpha_1^2  \\ 
\alpha_4   & \!\!\! =      B \, \Sigma \, \beta_3  + {1\over4} \, B_2 \, \beta_2 
+  {1\over6} \, B \, D_2 \, \Sigma \, \beta_1   
-   {1\over6} \, A \, B \, \beta_2  \\ & 
-  {1\over6} \, B \, E \, \alpha_1 \,  \alpha_2   -  {1\over6} \, B \, E \, \alpha_2  \,  \alpha_1 
- {1\over6} \, B \, \Sigma \, \beta_1  \,   \alpha_1^2  \, . 
\end {array} \right. \monend

\monitem Proof of Proposition 10.  

First observe that $ \, \dd \Phi . \, \xi = E \, \xi \, $ for any test vector $ \, \xi $. 
Due to the relation (\ref{Gamma-1}), we have 

\smallskip \noindent 
$ \Gamma_1 = A \, W + B \, \Phi(W) = A \, W + B \, E \, W = ( A + B \, E ) \, W \, $ 
and $ \, \alpha_1 =  A  + B \, E $. 
Due to (\ref{Phi-1}), we have moreover

\smallskip \noindent 
$ \Psi_1 = \dd \Phi(W) . \, \Gamma_1 (W) - \big( C \, W + D \, \Phi(W) \big) $

\smallskip \quad $ \,\, 
= E \, \alpha_1 \, W - C \, W - D \, E \,  W  $

\smallskip \quad $ \,\, 
= ( E \, \alpha_1 - C - D \, E ) \, W $ 

\smallskip \noindent 
and the second relation of (\ref{formules-lineaires}) is established. Due to (\ref{Gamma-2}),

\smallskip \noindent $
\Gamma_2 (W)  = B \, \Sigma \,\, \Psi_1  (W) $

\smallskip \qquad $ \quad  \, 
= B \, \Sigma \,\, \beta_1 \, W \, $. 

\smallskip \noindent 
Then  
$ \, \alpha_2 =  B \, \Sigma \,\, \beta_1 \, $ and the third line of (\ref{formules-lineaires}) is proven.
 We have now from (\ref{Phi-2}), 

\smallskip \noindent $
\Psi_2  (W)  = \Sigma \,  \dd \Psi_1 (W) . \Gamma_1 (W) + \dd \Phi(W) .  \Gamma_2 (W) - D \, \Sigma \,  \Psi_1 (W)  $

\smallskip \qquad $ \quad  \,\, =
 \Sigma \,  \beta_1 \, \alpha_1 \, W  + E \, \alpha_2 \, W - D \, \Sigma \,  \beta_1 \, W $  

\smallskip \qquad $ \quad  \,\, =
\big( \Sigma \,  \beta_1 \, \alpha_1   + E \, \alpha_2  - D \, \Sigma \,  \beta_1 \big) \, W $  

\smallskip \noindent 
and $ \, \beta_2 = \Sigma \,  \beta_1 \, \alpha_1   + E \, \alpha_2  - D \, \Sigma \,  \beta_1 \, $
as suggested in the fourth line of(\ref{formules-lineaires}).
From the relation (\ref{Gamma-3}), we have now

\smallskip \noindent $
\Gamma_3 (W)  = B \, \Sigma \, \Psi_2 (W) + {1\over12} \, B_2 \, \Psi_1 (W)
- {1\over6} \,B \, \dd \Psi_1 (W) .  \Gamma_1(W)  $ 

\smallskip \qquad $ \quad  \, =  
B \, \Sigma \, \beta_2 \, W  + {1\over12} \, B_2 \, \beta_1 \, W  - {1\over6} \, B \, \beta_1 .  \alpha_1 \, W  $ 

\smallskip \qquad $ \quad  \, =  
\big( B \, \Sigma \, \beta_2   + {1\over12} \, B_2 \, \beta_1 - {1\over6} \, B \, \beta_1 .  \alpha_1 \big) \, W  \, . $ 

\smallskip \noindent 
Then $ \, \alpha_3 = B \, \Sigma \, \beta_2   + {1\over12} \, B_2 \, \beta_1 - {1\over6} \, B \, \beta_1 .  \alpha_1 \, $
as proposed in the fifth line of (\ref{formules-lineaires}). 
Observe now that due to (\ref{notation-d2-psi-gamma1}), we have 

\smallskip \noindent $
\partial^2 \Psi_1 . \,  \Gamma_1 (W) =  \dd^2 \Psi_1 (W) . (\Gamma_1 ,\, \Gamma_1) 
+  \dd \Psi_1 (W) . \dd  \Gamma_1 (W)  .  \, \Gamma_1 $

\smallskip \qquad $ \quad  \, =  \dd \Psi_1 (W) . \dd  \Gamma_1 (W)  .  \, \Gamma_1 $ 
\hfill because $ \, \Psi_1 \, $ is linear 

\smallskip \qquad $ \quad  \, =  \beta_1 \, \alpha_1^2 \, W \, .$

\smallskip \noindent 
Then due to (\ref{Phi-3}), we have in this linear case, 

\smallskip \noindent $
\Psi_3  (W)  = \dd \Phi(W) .  \Gamma_3 (W)  +  \Sigma \,\,  \dd \Psi_1 (W) . \Gamma_2 (W) + \Sigma \,\, \dd \Psi_2 (W) . \Gamma_1 (W) - D \, \Sigma \,  \Psi_2 (W)$ 

\smallskip \qquad $ \qquad  \quad    
+ {1\over6} \, D \, \dd \Psi_1 (W) . \Gamma_1 (W) - {1\over12} \, D_2 \, \Psi_1 (W)
- {1\over12} \, \partial^2  \Psi_1(W) .  \Gamma_1 (W) $ 

\smallskip \qquad $ \quad  \,   =  
E \, \alpha_3 \, W + \Sigma \, \beta_1 \, \alpha_2 \, W 
 + \Sigma \, \beta_2 \, \alpha_1 \, W  - D \, \Sigma \,  \beta_2 \, W  
+ {1\over6} \, D \, \beta_1 \, \alpha_1 \, W $

\smallskip \qquad $ \qquad  \quad    - {1\over12} \, D_2 \, \beta_1 \, W - {1\over12} \, \beta_1 \, \alpha_1^2 \, W  $ 

\smallskip \qquad $ \quad  \,   =  
\big( E \, \alpha_3  + \Sigma \, \beta_1 \, \alpha_2  + \Sigma \, \beta_2 \, \alpha_1 - D \, \Sigma \,  \beta_2 
+ {1\over6} \, D \, \beta_1 \, \alpha_1   - {1\over12} \, D_2 \, \beta_1  - {1\over12} \, \beta_1 \, \alpha_1^2 \big) \, W $

\smallskip \noindent 
and  $ \, \beta_3 =  E \, \alpha_3  + \Sigma \, \beta_1 \, \alpha_2  + \Sigma \, \beta_2 \, \alpha_1 - D \, \Sigma \,  \beta_2 
+ {1\over6} \, D \, \beta_1 \, \alpha_1   - {1\over12} \, D_2 \, \beta_1  - {1\over12} \, \beta_1 \, \alpha_1^2 $.
We have finally, due to  (\ref{Gamma-4}), 

\smallskip \noindent $ 
\Gamma_4 (W)  = B \, \Sigma \, \Psi_3 (W) + {1\over4} \, B_2 \, \Psi_2 (W) + {1\over6} \,B \, D_2 \, \Sigma \, \Psi_1 - {1\over6} \, A \, B \, \Psi_2 (W)
- {1\over6} \, B \, \dd \gamma_1 (W) . \Gamma_2(W) $

\smallskip \qquad $ \qquad  \quad  
- {1\over6} \, B \, \dd \gamma_2 (W) . \Gamma_1 (W) - {1\over6} \, B \, \Sigma \,  \partial^2  \Psi_1(W) .  \Gamma_1 (W) \, $

\smallskip \qquad $ \,\,\,  =  
 B \, \Sigma \, \beta_3 \, W  + {1\over4} \, B_2 \, \beta_2 \, W + {1\over6} \,B \, D_2 \, \Sigma \, \beta_1 \, W - {1\over6} \, A \, B \, \beta_2 \, W 
- {1\over6} \, B \, E \, \alpha_2 \, W  $ 

\smallskip \qquad $ \qquad  \quad  
- {1\over6} \, B \, E \, \alpha_1 \, W - {1\over6} \, B \, \Sigma \,  \beta_1 \, \alpha_1^2  \, W \, $

\smallskip \qquad $ \,\,\,    =  
\big( B \, \Sigma \, \beta_3  + {1\over4} \, B_2 \, \beta_2  + {1\over6} \,B \, D_2 \, \Sigma \, \beta_1 - {1\over6} \, A \, B \, \beta_2  
- {1\over6} \, B \, E \, \alpha_2 \, \alpha_1 - {1\over6} \, B \, E \, \alpha_1 \, \alpha_2  
- {1\over6} \, B \, \Sigma \,  \beta_1 \, \alpha_1^2 \big) \, W  $

\smallskip \noindent 
and the proposition is established. \hfill $ \square $ 

\bigskip \bigskip    \noindent {\bf \large    8) \quad  Conclusion }    

\noindent
In this  
contribution, 
we have extended the Taylor expansion method  of a multiple relaxation times lattice Boltzmann scheme
up to fourth order accuracy. 
With this expansion, nonlinear partial differential equations for the conserved variables
and  various differential expressions for the nonconserved moments emerge naturally.
This expansion has been validated with previous works on the subject.

\smallskip
Observe that the main result is implemented  with this framework in the ``pyLBM''' software \cite{GG17}. 
The next step is the introduction of source terms, in the spirit proposed in \cite{DLT14}, 
and other geometries like triangles, in a way suggested in \cite{DL13}. 
The adaptation of the  Taylor expansion method to other variants of lattice
Boltzmann schemes as recentered schemes \cite{GGK06,DFG15b},
regularized lattice Boltzmann \cite{LC06,PH06,SHC06} or 
entropic lattice Boltzmann schemes  \cite{KGSB98}  are also into study.
After a first step in \cite{OBD18}, 
the link between the single scale Taylor expansion
and the multiple scale Chapman-Enskog is also a project for  the near future. 

%
\bigskip \bigskip   \noindent {\bf  \large  Acknowledgments }

\noindent 
The author thanks  Bruce Boghosian, Filipa Caetano, Lo\"ic Gouarin, Benjamin Graille, Pierre Lallemand
and Paulo Philippi for impressive disussions.
The scientific exchanges with the referee  greatly improved the matter of this contribution. 
All the formal computations have been implemented with SageMath \cite{sagemath}.

\bigskip \bigskip      \noindent {\bf  \large  References }   



\end{document}